\documentclass[preprint,12pt]{elsarticle}

\usepackage[width=16cm]{geometry}
\usepackage{amsmath}
\usepackage{amssymb}

\newcommand{\ud}{\textrm{d}}
\newcommand{\df}[2]{\frac{\partial #1}{\partial #2}}
\newcommand{\dd}[2]{\frac{\ud #1}{\ud #2}}
\newcommand{\re}{\mathbb{R}}
\newcommand{\norm}[1]{\|#1\|}

\begin{document}

\begin{frontmatter}

\title{Vertex-centroid finite volume scheme on tetrahedral grids for conservation laws}

\author{Praveen Chandrashekar}
\ead{praveen@math.tifrbng.res.in}
\author{Ashish Garg}
\ead{ashish\_garg31@yahoo.com}
\address{TIFR Center for Applicable Mathematics, Bangalore-560065, India}

\begin{abstract}
Vertex-centroid schemes are cell-centered finite volume schemes for conservation laws which make use of vertex values to construct high resolution schemes. The vertex values must be obtained through a consistent averaging (interpolation) procedure. A modified interpolation scheme is proposed which is better than existing schemes in giving positive weights in the interpolation formula. A simplified reconstruction scheme is also proposed which is also more accurate and efficient. For scalar conservation laws, we develop limited versions of the schemes which are stable in maximum norm by constructing suitable limiters. The schemes are  applied to compressible flows governed by the Euler equations of inviscid gas dynamics.
\end{abstract}

\begin{keyword}
Finite volume method, unstructured grids, reconstruction, maximum principle, compressible flows
\end{keyword}

\end{frontmatter}

\section{Introduction}

Finite volume methods on unstructured grids have become very useful in solving conservation laws, especially those arising in fluid dynamics and in which the computational domain is of a complex shape. This is mainly due to the maturity of grid generation tools, especially triangular and tetrahedral grids. Also, the local conservation property of finite volume methods allows them to automatically compute discontinuous solutions. Finite volume methods can be classified based on the location where the solution values are stored~\cite{barth94aspects,mavriplis2007}. In {\em cell-centered} methods, the solution is assumed to be stored at the centroid of the cells, while in {\em vertex-centered} methods, the solution is stored at the vertices of the grid. In the latter case, a cell has to be constructed around each vertex for which there are several possible alternatives. The most commonly used method is the dual grid scheme in which the cell is obtained by joining the cell centroids to the face and edge centroids. Strictly speaking, the unknowns in a finite volume method are cell averages, which are not associated with any particular spatial location in the cell. However, the cell average value gives a second order accurate approximation to the solution at the cell centroid. For constructing second order accurate schemes, the cell average value can be taken to be the point value at the cell center, which is the approach taken in the present work.

The construction of second order accurate schemes leads to several possible alternatives. Since our work uses cell centered schemes, we restrict the discussion to this case. The finite volume method updates the cell average value of the solution, and the detailed variation of the solution inside the cell is not known. To construct higher order schemes, the solution inside each cell must be {\em reconstructed} by using the cell average values of the neigbouring cells. The common approach is to reconstruct the solution variation in each cell by a linear polynomial by making use of the cell average values in the current cell and some of its neighbours, which forms the stencil of the reconstruction. This can be achieved by first approximating the derivatives of the solution in the cell by using either a Green-Gauss procedure or a least squares procedure~\cite{barth94aspects}. The reconstructed solution may be limited in order to maintain oscillation-free solutions. While in one dimensional problems, the TVD criterion~\cite{harten1984} can be used, for multi-dimensional problems on unstructured grids, the usual criterion is to ensure that the reconstructed solution lies within the minimum and maximum of the neigbouring cell average values. The two states obtained on either side of a cell face by the reconstructed solution are used to compute the flux across the face using some numerical flux function.

The vertex-centroid scheme~\cite{frink1994,frink1998} is a finite volume scheme for conservation laws on triangles and tetrahedral grids. The basic unknowns are the cell center values but it also makes use of vertex values in order to construct the second order scheme.  The authors in~\cite{frink1994,frink1998} obtain a formula for the reconstructed value at the center of a cell face by making use of the geometrical properties of triangles/tetrahedra, and using the cell center and vertex values of the solution.  There is no need to explicitly compute the derivatives of the solution in each cell as in the usual cell-centered finite volume schemes. The vertex values are obtained through an interpolation/averaging procedure that is exact for linear polynomials. This procedure has been refered to as {\em pseudo-laplacian averaging} in the literature. The interpolation formula for each vertex is a linear combination of the cell-center values adjacent to that vertex and it is desirable from a stability point of view that the weights in this formula should be positive. Applications to compressible transonic flows with weak shocks have been given in~\cite{frink1994,frink1998} and the vertex-centroid scheme was able to compute stable and accurate solutions even without a limiter, demonstrating good stability properties of the scheme. For stronger shocks, the authors recommend using a minmod limiter though a stability analysis was not shown.

A vertex-centroid scheme which is similar to the schemes in~\cite{frink1994,frink1998} has been used in~\cite{jameson-iccfd2000} for two dimensional compressible flows, but it uses an area averaging procedure which is not second order accurate on general meshes. However the authors use very regular meshes and find that the solutions have extremely low levels of numerical dissipation. They also comment that the vertex-centroid scheme simplified the inter-mesh transfers of an unstructured mesh multi-grid scheme, in both the fine-to-coarse and coarse-to-fine directions.

In this paper, we analyze and improve the vertex-centroid schemes in three respects. We present a modified averaging procedure which reduces the number of negative weights. With the new scheme the matrices arising in the determination of the weights have better behaviour. Secondly, we analyze the accuracy of Frink reconstruction procedure and also propose a simplied procedure which is theoretically more accurate. Thirdly, we develop limited versions of the schemes which can be shown to be stable in maximum norm when applied to scalar conservation laws. Under a CFL condition, the solution at the new time level in each cell is shown to be bounded between the solution at the neigbouring cells and the average values at the vertices of the cell. This is an important design criterion for numerical schemes solving hyperbolic conservation laws since they can develop discontinuous solutions, even starting from smooth initial data. The developed schemes are tested on several problems governed by the inviscid Euler equations modeling compressible flows including those with discontinuous solutions. The numerical solutions are compared with exact solutions or with experimental data to demonstrate the performance of the schemes.

The rest of the paper is organized as follows. The basic formulation of the vertex-centroid scheme is explained. The vertex interpolation schemes are discussed and a new interpolation scheme is proposed. We then study the error in the reconstruction schemes of Frink and a simpler reconstruction scheme proposed here which we refer to as {\em upwind reconstruction}. Limited versions of both the reconstruction schemes are described and their maximum stability is analyzed for a scalar conservation law in two and three dimensions. Finally, the schemes are applied to problems governed by Euler equations of gas dynamics and their performance in computing discontinuous solutions is demonstrated.
\section{Vertex-centroid scheme}

Consider a system of conservation laws in $d$ space dimensions which can be written as
\begin{equation}
\df{U}{t} + \sum_{i=1}^d \df{F_i}{x_i} = 0
\end{equation}
Here $U$ is the vector of conserved variables and $F_i$, $i=1, \ldots, d$ are the Cartesian components of the flux vector and $d$ is the number of spatial dimensions. We solve this equation numerically on a triangulation made of triangles in 2-D and tetrahedra in 3-D. We assume that the cells and vertices are numbered in some way; the $i$'th cell is denoted by $C_i$. The face between the $i$'th cell and the $j$'th cell is denoted $S_{ij}$. The set of indices of the cells sharing a face with $C_i$ is denoted by $N(i)$. Then the semi-discrete finite volume scheme for cell $C_i$ is
\begin{equation}
|C_i| \dd{U_i}{t} + \sum_{j \in N(i)} H(U^+_{ij},U^-_{ij},n_{ij})=0
\label{eq:semid}
\end{equation}
where $n_{ij} \in \re^d$ is the normal to $S_{ij}$ and pointing into cell $C_j$ with $|n_{ij}|=|S_{ij}|$, $|\cdot|$ denotes the measure of the corresponding geometric objects and $H$ is a numerical flux function.   Note that we have used a mid-point quadrature rule to approximate the flux integral on the cell faces which is sufficient to obtain second order accuracy. The numerical flux function must satisfy the following conditions:
\begin{itemize}
\item Consistency
\begin{equation*}
H(U,U,n)=\sum_{i=1}^d F_i(U) n_i
\end{equation*}
\item Conservation
\begin{equation*}
H(U_1,U_2,n) + H(U_2,U_1,-n)=0
\end{equation*}
\end{itemize}
In the scalar conservation law case, the numerical flux function is a monotone;  in this case, the numerical flux $H(a,b,n)$ is an increasing function of the first argument and a decreasing function of the second argument.

\begin{figure}
\begin{center}
\includegraphics[width=0.5\textwidth]{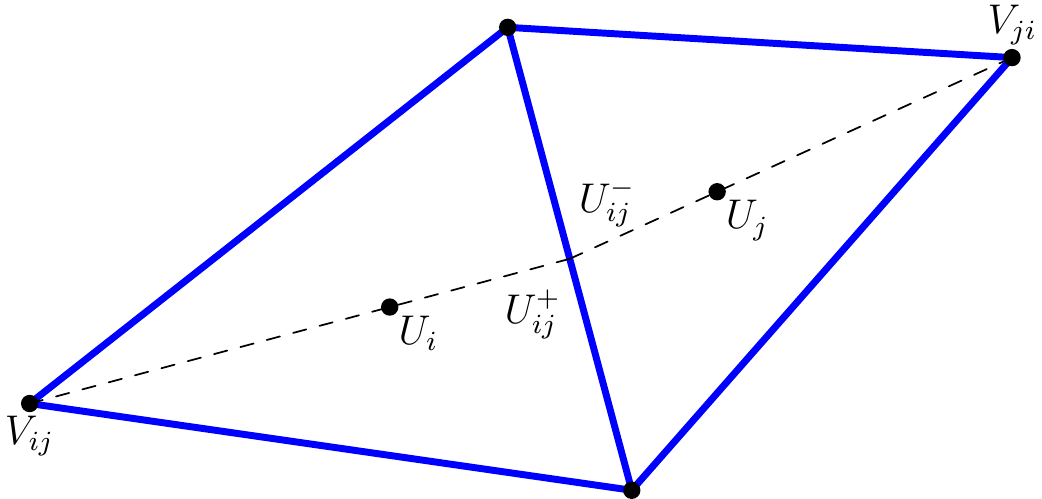}
\caption{Reconstruction in vertex-centroid scheme}
\label{fig:r2}
\end{center}
\end{figure}

The values $U_{ij}^+$ and $U_{ij}^-$ are the two states at the center of face $S_{ij}$ obtained by some reconstruction process in the cells $C_i$ and $C_j$ respectively, see figure~(\ref{fig:r2}). In any triangle/tetrahedron, the line joining any vertex to the center of the opposite face pasess through the centroid of the cell. The ratio of the distance between any vertex and the cell centroid, to the distance between the cell centroid and the center of the opposite face is 2:1 for a triangle and 3:1 for a tetrahedron; these properties are used in the reconstruction schemes. Let $V_{ij}$ be the value at the vertex of cell $C_i$ which is opposite to the face $S_{ij}$ and $W_{ij}$ denote the arithmetic average of the vertex values on the face $S_{ij}$. Note that $W_{ij}=W_{ji}$ since both the values correspond to the same face, but $V_{ij} \ne V_{ji}$ since they refer to different vertices. Then the reconstructed values on either side of a face $S_{ij}$ is given by Frink~\cite{frink1994} as
\begin{equation}
U^+_{ij} = U_i + \alpha_d (W_{ij} - V_{ij}), \qquad U^-_{ij} = U_j + \alpha_d (W_{ij} - V_{ji})
\label{eq:frink}
\end{equation}
where the subscript $d$ refers to the spatial dimension; for 2-D we have $\alpha_2 = 1/3$ while in 3-D, we have $\alpha_3 = 1/4$. We will refer to this as the {\em Frink scheme}.

An alternate reconstruction which does not make use of face averages is given by
\begin{equation}
U^+_{ij} = U_i + \beta_d (U_i - V_{ij}), \qquad U^-_{ij} = U_j + \beta_d (U_j - V_{ji})
\label{eq:new}
\end{equation}
where $\beta_2 = 1/2$ and $\beta_3=1/3$. We will refer to this as the {\em upwind scheme}. This scheme is a simple consequence of extrapolating the solution at the cell center and vertex to the center of the opposite face. Since this scheme uses states on only one side of the face for the reconstruction, we refer to this as the upwind reconstruction scheme. Computationally, the scheme given by equation~(\ref{eq:new}) is more simple and efficient compared to (\ref{eq:frink}) since it does not require the face average values $W_{ij}$. Later we will see that the upwind scheme is also advantageous when limiters are used, since it leads to less restrictions on the reconstruction to be stable in maximum norm.
\section{Vertex interpolation scheme}

The basic finite volume scheme used here is based on cell center values which are updated by the finite volume scheme. In this case, we do not know the solution at the vertices of the grid, which has to be approximated through some consistent interpolation procedure. To discuss the interpolation scheme, we take a generic vertex and assume without loss of generality, that it is located at the origin. We then find all the cells which contain this vertex; let the centroids of these cells have coordinates $(x_j,y_j,z_j)$, $j=1,\ldots,n$ and let $r_j$ be the Euclidean distance from the $j$'th centroid to the vertex. In the literature, the following interpolation procedures have been used, an area/volume averaging procedure and an inverse distance averaging procedure, given by
\begin{displaymath}
V = \frac{\sum_j |C_j| U_j}{\sum_j |C_j|} \qquad \textrm{or} \qquad V = \frac{\sum_j \frac{1}{r_j} U_j}{\sum_j \frac{1}{r_j}}
\end{displaymath}
where the summations are over all the $n$ cells containing the vertex. The inverse distance averaging can be viewed as a local Shepard interpolation~\cite{Shepard:1968:TIF:800186.810616}. On general meshes, both the above procedures are only first order accurate, i.e., they are exact only for constant functions. On uniform and isotropic grids, they become second order accurate. Note that for the finite volume scheme to be second order accurate the interpolation scheme must be atleast second order accurate. In~\cite{jameson-iccfd2000}, the area averaging procedure has been used together with high quality meshes leading to very accurate solutions. While these methods may not be accurate in general, they are very robust since they are convex combinations of the cell center values and hence the interpolated value is bounded between the minimum and maximum of the cell center values. This is a useful property while solving hyperbolic conservation laws which can have discontinuous solutions.

\subsection{Laplacian averaging}

In order to construct an interpolation procedure which is exact for linear polynomials, Frink~\cite{frink1994} extended a method for 2-D Navier-Stokes solutions given in \cite{holmes1989} to three dimensions. In this approach, we first assume an interpolation formula of the type
\begin{equation}
V = \frac{\sum_j w_j U_j}{\sum_j w_j}
\label{eq:plapavg}
\end{equation}
The necessary conditions for this formula to be exact for linear polynomials are
\begin{equation}
\sum_j w_j x_j =0, \quad \sum_j w_j y_j = 0, \quad \sum_j w_j z_j = 0
\label{eq:avgcon}
\end{equation}
In order to determine the weights $w_j$, the following minimization problem is solved with respect to the weights together with the constraints given by equations~(\ref{eq:avgcon})
\begin{equation}
\min_{\{w_j\}} \frac{1}{2} \sum_j (w_j - 1)^2
\label{eq:wtfun}
\end{equation}
This approach tries to keep the weights as close to unity as possible while still satisfying the consistency conditions. In order to solve the constrained minimization problem, we introduce Lagrange multipliers $\lambda_x, \lambda_y, \lambda_z$; the constrained minimization problem is converted into an unconstrained problem in which we have to now minimize the function
\begin{equation*}
 \frac{1}{2} \sum_j (w_j - 1)^2 + \lambda_x \sum_j w_j x_j + \lambda_y \sum_j w_j y_j + \lambda_z \sum_j w_j z_j
\end{equation*}
with respect to the weights $\{w_j\}$ and the Lagrange multipliers $\lambda_x, \lambda_y, \lambda_z$ as independent variables. The solution is obtained by setting the derivative of the above function with respect to each of these variables to zero; the weights are then given by
\begin{displaymath}
w_j = 1 + \lambda_x x_j + \lambda_y y_j + \lambda_z z_j
\end{displaymath}
where the Lagrange multipliers are obtained by solving the following system of equations
\begin{equation}
\begin{bmatrix}
\sum x_j^2 & \sum x_j y_j & \sum x_j z_j \\
\sum x_j y_j & \sum y_j^2 & \sum y_j z_j \\
\sum x_j z_j & \sum y_j z_j & \sum z_j^2
\end{bmatrix}
\begin{bmatrix}
\lambda_x \\ \lambda_y \\ \lambda_z
\end{bmatrix} = -
\begin{bmatrix}
\sum x_j \\ \sum y_j \\ \sum z_j
\end{bmatrix}
\label{eq:lapmat}
\end{equation}
In \cite{frink1994} the above set of equations is solved explictly using Cramers rule.
\subsection{Consistent Shepard interpolation}

We propose a modification of the Shepard interpolation procedure which makes it exact for linear polynomials. In this case, we begin by assuming an interpolation formula of the form
\begin{equation}
V = \frac{\sum_j \frac{w_j}{r_j} U_j}{\sum_j \frac{w_j}{r_j}}
\label{eq:avgshep}
\end{equation}
where the weights should now satisfy the consistency conditions
\begin{equation}
\sum_j w_j \frac{x_j}{r_j} =0, \quad \sum_j w_j \frac{y_j}{r_j} = 0, \quad \sum_j w_j \frac{z_j}{r_j} = 0
\label{eq:avgcon2}
\end{equation}
In order to determine the weights $w_j$, the minimization probem (\ref{eq:wtfun}) is solved with respect to the weights together with the constraints given by equations~(\ref{eq:avgcon2}).
In order to solve the constrained minimization, we again introduce Lagrange multipliers $\lambda_x, \lambda_y, \lambda_z$ in terms of which the solution can be written as
\begin{displaymath}
w_j = 1 + \lambda_x \frac{x_j}{r_j} + \lambda_y \frac{y_j}{r_j} + \lambda_z \frac{z_j}{r_j}
\end{displaymath}
where the Lagrange multipliers are obtained by solving the following system of equations
\begin{equation}
\begin{bmatrix}
\sum \left(\frac{x_j}{r_j}\right)^2 & \sum \left(\frac{x_j}{r_j}\right) \left(\frac{y_j}{r_j}\right) & \sum \left(\frac{x_j}{r_j}\right) \left(\frac{z_j}{r_j}\right) \\
\sum \left(\frac{x_j}{r_j}\right) \left(\frac{y_j}{r_j}\right) & \sum \left(\frac{y_j}{r_j}\right)^2 & \sum \left(\frac{y_j}{r_j}\right)  \left(\frac{z_j}{r_j}\right) \\
\sum \left(\frac{x_j}{r_j}\right) \left(\frac{z_j}{r_j}\right) & \sum \left(\frac{y_j}{r_j}\right) \left(\frac{z_j}{r_j}\right) & \sum \left(\frac{z_j}{r_j}\right)^2
\end{bmatrix}
\begin{bmatrix}
\lambda_x \\ \lambda_y \\ \lambda_z
\end{bmatrix} = -
\begin{bmatrix}
\sum \left(\frac{x_j}{r_j}\right) \\ \sum \left(\frac{y_j}{r_j}\right) \\ \sum \left(\frac{z_j}{r_j}\right)
\end{bmatrix}
\label{eq:shepmat}
\end{equation}
We solve the above matrix equation explicitly using Cramer's rule.
\paragraph{Remarks}

A first observation in comparing the above two methods is that all the formulae in the modified Shepard interpolation are independent of the spatial scales. If $h$ is a measure of the local grid size, then the determinant of the matrix in equation~(\ref{eq:lapmat}) is $O(h^6)$ and the numerical value of this determinant can become very small. For example, if $h=10^{-2}$, then the determinant could be $O(10^{-12})$. In the examples, we find that the determinant can reach the levels of machine precision even for inviscid grids. Due to this reason, the solution of equation~(\ref{eq:lapmat}) will also suffer from round-off errors. On the other hand, the determinant of the matrix in equation~(\ref{eq:shepmat}) is independent of $h$ and the determinant is found to be well behaved, as given in the later examples.

Both of the above methods do not guarantee that the weights $w_j$ will be positive. A necessary condition for equations~(\ref{eq:avgcon}) or (\ref{eq:avgcon2}) to be satisfied with positive weights is that the cell centers must be distributed on all sides of the vertex. This means that there must be atleast one cell with $x_j < 0$ and atleast one cell with $x_j >0$, with similar conditions for the other two coordinates. Also, these conditions must be satisfied for any orientation of the coordinate axes. In equations~(\ref{eq:avgcon}) the magnitude of the different terms depends on the spatial distribution of the cell centers around the vertex while in (\ref{eq:avgcon2}) the magnitudes depend only on the angular distribution of the cell centers. We conjecture that due to this reason, the modified Shepard interpolation scheme will be able to preserve the positivity of the weights better than the Laplacian averaging scheme. We indeed find this to be the case in our numerical tests given in the later sections. In many cases, the modified scheme gives all positive weights while the Laplacian scheme yields some negative weights. In other cases, the number of vertices with negative weights is considerably reduced with the modified scheme.

Due to the presence of the inverse distance factor in equation~(\ref{eq:avgshep}), the modified Shepard interpolation scheme gives more weightage to the cells which are closer to the given vertex. This might be beneficial for anisotropic grids which would arise if we perform anisotropic  grid adaptation, or in the case of boundary layer grids. However in the present work, we have only considered inviscid problems for which the grids are not highly anistropic.
\section{Accuracy of reconstruction schemes}

The derivation of the reconstruction scheme (\ref{eq:frink}) makes use of the geometrical properties of triangles and tetrahedra, leading to second order accuracy. This scheme makes use of all the vertex values including the cell center value. A simpler scheme which is also second order accurate is given by equation~(\ref{eq:new}). This is based on a one dimensional extrapolation of the vertex and cell center values to the mid-point of the opposite face. In this section we study the accuracy of these two schemes by assuming that the vertex values are exact.

Assume that $U$ is twice continuously differentiable function of the space coordinates. We denote by $D^2U$ the symmetric matrix of second derivatives. In any cell $C$, we have
\begin{equation}
\norm{D^2U} \le K
\end{equation}
for some positive constant $K$, where the matrix norm $\|\cdot\|$ is induced by the $l_2$ norm for vectors~\cite{kreyszig}.
\subsection{2-D reconstruction}

\begin{figure}
\begin{center}
\includegraphics[width=0.4\textwidth]{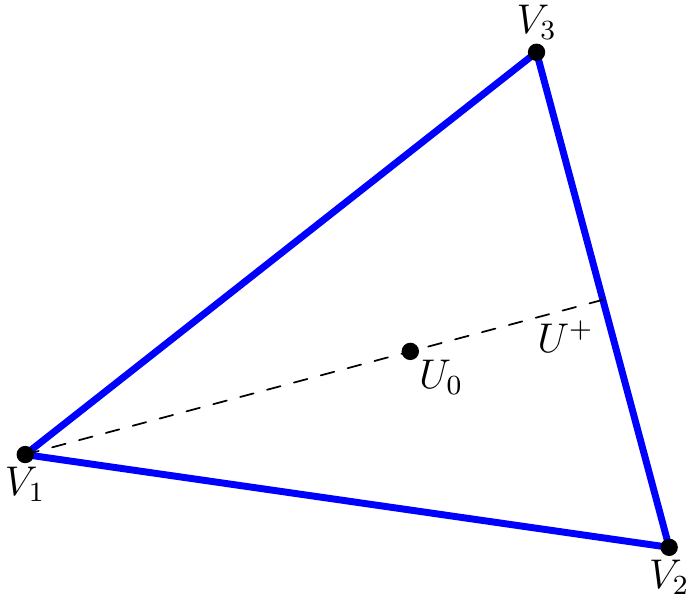}
\includegraphics[width=0.35\textwidth]{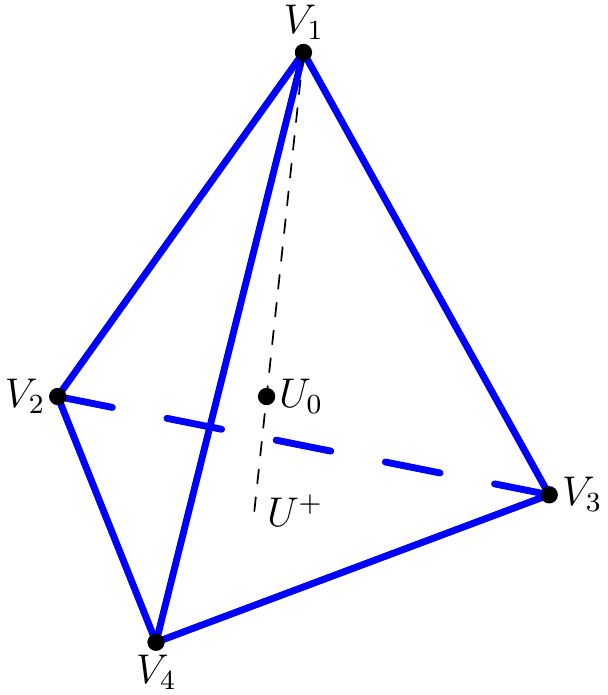}
\caption{2-D and 3-D reconstruction in vertex-centroid scheme}
\label{fig:r3}
\end{center}
\end{figure}

Consider a triangle $C$ whose vertex values are $V_1,V_2,V_3$ and whose cell center value is $U_0$, see figure~(\ref{fig:r3}). Without loss of generality, assume that the cell center is at the origin; then the position vector of the vertices $r_1, r_2, r_3 \in \re^2$ are such that $r_1 + r_2 + r_3 = 0$. Let us look at the reconstructed value on the face formed by vertices 2 and 3. The position vector of the center of this face is $\frac{1}{2}(r_2 + r_3)=-\frac{1}{2}r_1$. The exact value at the face center is given by a Taylor series with remainder~\cite{rudin1} as
\begin{equation}
U_e = U_0 - \frac{1}{2} r_1^T G_0 + \frac{1}{8}r_1^T H_{01} r_1
\label{eq:ueexact2}
\end{equation}
where we denote $G_0 = \nabla U_0$ and $H_{01}=D^2U(-s r_1/2)$ for some $s \in [0,1]$. Using another Taylor series with remainder term, we can write the vertex values in terms of the cell center values as
\begin{equation}
V_j = U_0 + r_j^T G_0 + \frac{1}{2} r_j^T H_j r_j
\label{eq:v2d}
\end{equation}
where  $H_j = D^2 U(s_j r_j) \in \re^{2 \times 2}$ for some $s_j \in [0,1]$ and $j=1,2,3$.
\subsubsection{Frink scheme}

The reconstructed value on the face formed by vertices 2, 3 is according to equation~(\ref{eq:frink})
\begin{equation}
U^+ = U_0 + \frac{1}{3} \left[ \frac{1}{2}(V_2 + V_3) - V_1 \right]
\label{eq:frink2}
\end{equation}
which upon using equation~(\ref{eq:v2d}) leads to
\begin{equation*}
U^+ = U_0 - \frac{1}{2} r_1^T G_0 + \frac{1}{12} r_2^T H_2 r_2 + \frac{1}{12} r_3^T H_3 r_3 - \frac{1}{6} r_1^T H_1 r_1
\end{equation*}
Using equation~(\ref{eq:ueexact2}, the error in the interpolated value can be bounded as
\begin{equation*}
|U^+-U_e| \le \frac{7}{24} K \norm{r_1}^2 + \frac{1}{12} K \norm{r_2}^2 + \frac{1}{12} K \norm{r_3}^3 \le \frac{11}{24} K h^2
\end{equation*}
where the diameter of the cell is defined as
\begin{equation*}
h = \max_{1\le j \le 3} \norm{r_j}
\end{equation*}
\subsubsection{Upwind scheme}

The reconstructed value is given by equation~(\ref{eq:new}) as
\begin{equation}
U^+ = U_0 + \frac{1}{2} \left[ U_0 - V_1 \right]
\label{eq:up2}
\end{equation}
which upon using equation~(\ref{eq:v2d}) leads to
\begin{equation*}
U^+ = U_0 - \frac{1}{2} r_1^T G_0 - \frac{1}{4} r_1^T H_1 r_1
\end{equation*}
The error in this approximation can be bounded as
\begin{equation*}
|U^+ - U_e| \le \frac{3}{8} K \norm{r_1}^2 \le \frac{3}{8} K h^2
\end{equation*}
We see that the upwind scheme given by equation~(\ref{eq:up2}) has smaller error constant than the Frink scheme given by equation~(\ref{eq:frink2}).

\subsection{3-D reconstruction}

Consider a tetrahedron $C$ whose vertex values are $V_1,V_2,V_3,V_4$ and whose cell center value is $U_0$, see figure~(\ref{fig:r3}). Without loss of generality, assume that the cell center is at the origin; then the  position vectors of the vertices $r_1, r_2, r_3, r_4 \in \re^3$ are such that $r_1 + r_2 + r_3 + r_4= 0$. Let us look at the reconstructed value on the face formed by vertices 2, 3 and 4. The position vector of the center of this face is $\frac{1}{3}(r_2 + r_3 + r_4)=-\frac{1}{3}r_1$. The exact value at the face center is given by a Taylor expansion with remainder term as
\begin{equation}
U_e = U_0 - \frac{1}{3} r_1^T G_0 + \frac{1}{18}r_1^T H_{01} r_1
\label{eq:ueexact3}
\end{equation}
where we denote $G_0=\nabla U_0$ and $H_{01}=D^2U(-s r_1/3)$ for some $s \in [0,1]$. Using another Taylor series with remainder term, we can write the vertex values in terms of the cell center values as
\begin{equation}
V_j = U_0 + r_j^T G_0 + \frac{1}{2} r_j^T H_j r_j
\label{eq:v3d}
\end{equation}
where $H_j = D^2 U(s_j r_j) \in \re^{3 \times 3}$ for some $s_j \in [0,1]$ and $j=1,\ldots,4$.
\subsubsection{Frink scheme}

The reconstructed value on the face formed by vertices 2, 3 and 4 according to equation~(\ref{eq:frink}) is
\begin{equation}
U^+ = U_0 + \frac{1}{4} \left[ \frac{1}{3}(V_2 + V_3 + V_4) - V_1 \right]
\label{eq:frink3}
\end{equation}
which upon using equation~(\ref{eq:v3d}) leads to
\begin{equation*}
U^+ = U_0 - \frac{1}{3} r_1^T G_0 + \frac{1}{24} r_2^T H_2 r_2 + \frac{1}{24} r_3^T H_3 r_3 + \frac{1}{24} r_4^T H_4 r_4 - \frac{1}{8} r_1^T H_1 r_1
\end{equation*}
Using equation~(\ref{eq:ueexact3}), the error in the interpolated value can be bounded as
\begin{equation*}
|U^+-U_e| \le \frac{13}{72} K \norm{r_1}^2 + \frac{1}{24} K \norm{r_2}^2 + \frac{1}{24} K \norm{r_3}^3 + \frac{1}{24} K \norm{r_4}^2 \le \frac{11}{36} K h^2
\end{equation*}
where the cell diameter is defined as
\begin{equation}
h = \max_{1\le j \le 4} \norm{r_j}
\end{equation}
\subsubsection{Upwind scheme}

The reconstructed value is given by equation~(\ref{eq:new}) as
\begin{equation}
U^+ = U_0 + \frac{1}{3} \left[ U_0 - V_1 \right]
\label{eq:up3}
\end{equation}
which upon using equation~(\ref{eq:v3d}) leads to
\begin{equation*}
U^+ = U_0 - \frac{1}{3} r_1^T G_0 - \frac{1}{6} r_1^T H_1 r_1
\end{equation*}
Using equation~(\ref{eq:ueexact3}), the error in this approximation can be bounded as
\begin{equation*}
|U^+ - U_e| \le \frac{2}{9} K \norm{r_1}^2 \le \frac{2}{9} K h^2
\end{equation*}
We see that the upwind scheme given by equation~(\ref{eq:up3}) has smaller error constant than the Frink scheme given by equation~(\ref{eq:frink3}).

\paragraph{Remark} This analysis shows that the simpler upwind scheme is actually slightly more accurate compared to the Frink scheme. However both methods are second order accurate. The difference in the error constant does not translate into any appreciable difference in the numerical solutions obtained by the two methods, atleast for the test cases that we have studied. The purpose of this analysis was to show that there is no advantage in using the more elaborate scheme given by equation~(\ref{eq:frink}) and the more simpler scheme of equation~(\ref{eq:new}) is more accurate in theory. In the next section where we construct schemes stable in the maximum norm, we find that the upwind scheme is more beneficial since it requires less restrictions to satisfy the maximum stability.

\section{Maximum stability for scalar conservation law}

Solutions of nonlinear hyperbolic PDEs can develop discontinuous even if the initial condition is smooth. While first order accurate finite volume methods  on unstructured grids are capable of computing discontinuous solutions in a stable manner~\cite{BarJ89}, the higher order numerical schemes might produce spurious oscillations near discontinuities. The classic approach is to construct total variation diminishing schemes which will prevent the generation of oscillations. The TVD concept however is difficult to extend to multiple dimensions and unstructured grids. Moreover, there is a negative result by Goodman and LeVeque~\cite{goodman-leveque} that a TVD scheme can be at most first order accurate in more than one spatial dimension. Due to these reasons, one constructs schemes which are stable in the maximum norm. This prevents the generation of new extrema and hence oscillations in the solution will be avoided. The first order accurate finite volume schemes combined with a monotone flux are stable in the maximum norm. The basic idea to achieve this stability for higher order schemes is to reduce the accuracy of reconstruction to first order wherever the solution has a discontinuity and might lead to oscillations.

In the case of structured grids, there are well established conditions that can be used to check the positivity of the schemes. For one dimensional probems, we can write the second order semi-discrete scheme as
\begin{equation*}
h \dd{U_i}{t} = A_{i} (U_{i+1} - U_i) + B_{i} (U_{i-1} - U_i)
\end{equation*}
If all the coefficients $A$, $B$ are positive, then the scheme will be stable in the maximum norm under a CFL condition. In particular, if $U_i$ is a local maximum, then it will not increase in time, and similarly, if $U_i$ is a local minimum, it will not decrease in time. For unstructured grids, the cells are not arranged along grid lines which makes it difficult to write the scheme in the above form. However, for the vertex-centroid scheme, this can be achieved by using the vertex values along with the cell-center values, which is the key step in the present construction of limited schemes.

We propose to write the semi-discrete vertex-centroid finite volume scheme in the following form, which not only contains the difference of cell average values but also the difference between vertex value and cell average value.
\begin{equation}
|C_i| \dd{U_i}{t} = \sum_{j \in N(i)} \left[ A_{ij} (U_j - U_i) + B_{ij} (V_{ij} - U_i) \right]
\label{eq:vc}
\end{equation}
We discretize the above set of ordinary differential equations using Euler time stepping scheme which leads to the following fully discrete update equation
\begin{equation}
U^{n+1}_i = \left[1 - \frac{\Delta t}{|C_i|}\sum_{j \in N(i)} (A_{ij}+B_{ij}) \right] U^n_i + \frac{\Delta t}{|C_i|} \sum_{j \in N(i)} A_{ij} U_j + \frac{\Delta t}{|C_i|} \sum_{j \in N(i)} B_{ij} V_{ij}
\end{equation}
 If all the coefficients $A, B$ in the above equation are non-negative, then the resulting scheme satisfies a maximum principle under a CFL condition, i.e.,
\begin{equation}
\min_{j \in N(i)} \{ U^n_i, U^n_j, V^n_{ij} \} \le U^{n+1}_i \le \max_{j \in N(i)} \{U^n_i, U^n_j, V^n_{ij} \}
\end{equation}
 The vertex values are obtained by an interpolation formula of the type given in equation (\ref{eq:plapavg}) or (\ref{eq:avgshep}). If all the weights in the interpolation formula are positive, then the vertex values are bounded by the cell average values. The higher order accurate scheme must be constructed to ensure the positivity of the coefficients in equation (\ref{eq:vc}) so that spurious oscillations near shocks are prevented. If the second order reconstruction given by equation~(\ref{eq:frink}) or equation~(\ref{eq:new}) is written as
\begin{equation}
U^+_{ij}=U_i + \Delta U_{ij}, \qquad U^-_{ij}=U_j + \Delta U_{ji}
\end{equation}
then the corresponding limited reconstruction scheme is defined to be
\begin{equation}
U^+_{ij}=U_i + \theta_{ij} \Delta U_{ij}, \qquad U^-_{ij}=U_j + \theta_{ji} \Delta U_{ji}
\end{equation}
where the limiter $\theta \in [0,1]$ should be chosen to satisy the positivity conditions. When $\theta=1$, the scheme is second order accurate while if $\theta=0$, it becomes first order accurate.

Consider a numerical flux function which can be written as
\begin{equation}
H(a,b,n)=H^+(a,n) + H^-(b,n)
\end{equation}
with $H^+$ being a non-decreasing function and $H^-$ being a non-increasing function. Then from the consistency of the numerical flux function,
\begin{displaymath}
\sum_{j \in N(i)} H(U_i,U_i,n_{ij}) = 0
\end{displaymath}
we can separate the positive and negative fluxes and write equation~(\ref{eq:semid}) as
\begin{eqnarray*}
|C_i| \dd{U_i}{t} &=& -\sum_{j \in N(i)} [H(U^+_{ij},U^-_{ij},n_{ij})- H(U_i,U_i,n_{ij})] \\
&=&-\sum_{j \in N(i)} \underbrace{\frac{H^+(U^+_{ij},n_{ij})-H^+(U_i,n_{ij})}{U^+_{ij}-U_i}}_{P_{ij} \ge 0} \theta_{ij} \Delta U_{ij} \\
&& -\sum_{j \in N(i)} \underbrace{\frac{H^-(U^-_{ij},n_{ij})-H^-(U_i,n_{ij})}{U^-_{ij}-U_i}}_{Q_{ij} \le 0} (U_j + \theta_{ji}  \Delta U_{ji} - U_i)
\end{eqnarray*}
Due to the properties of the flux function, we know the sign of the terms defined as $P_{ij}$ and $Q_{ij}$. We next derive conditions on the limiter function $\theta$ for the two reconstruction schemes so that when they are written in the form of equation~(\ref{eq:vc}), they have non-negative coefficients.

\subsection{Frink scheme}

\begin{figure}
\begin{center}
\includegraphics[width=0.5\textwidth]{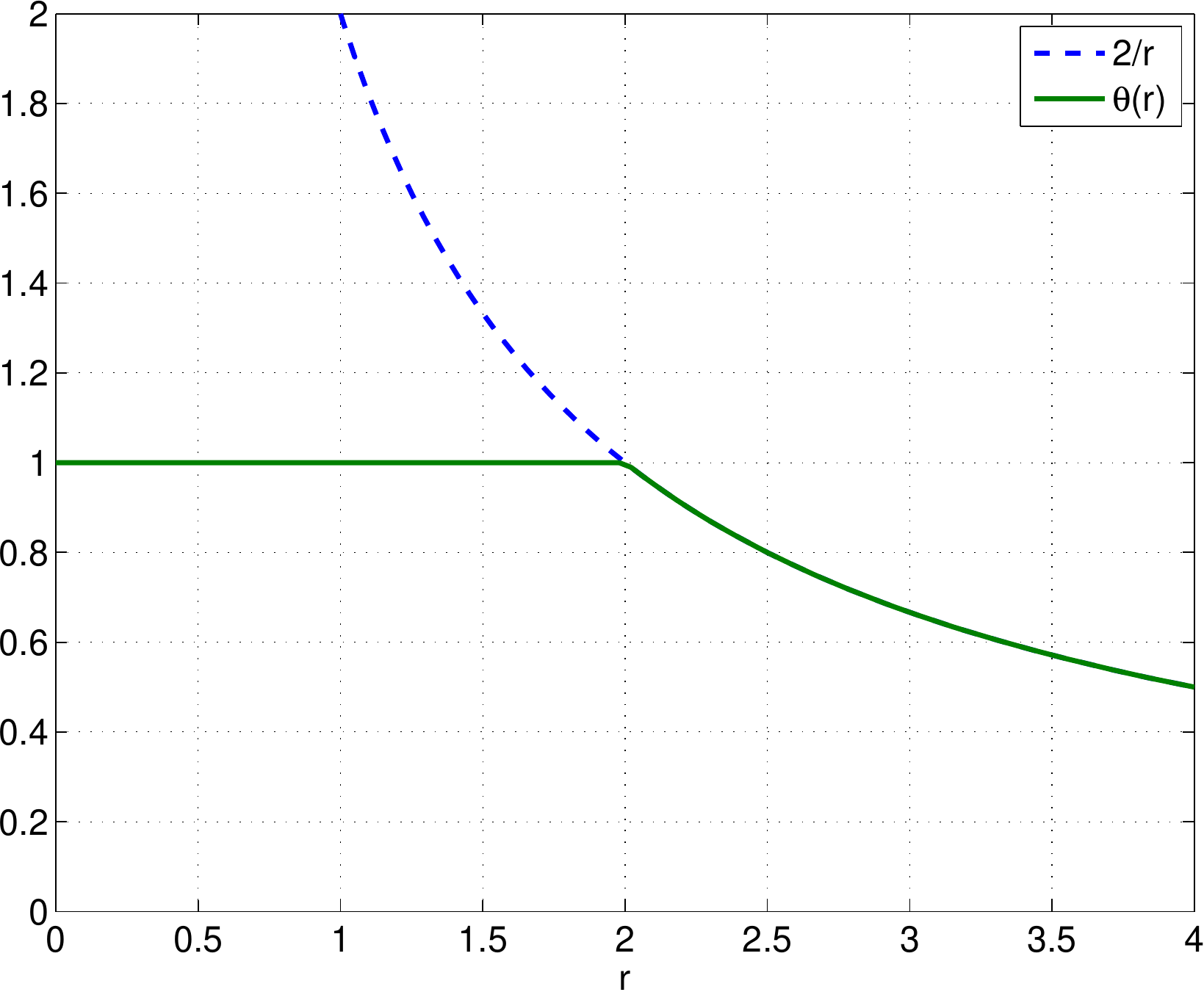}
\caption{Limiter function $\theta(r)$}
\label{fig:limiter}
\end{center}
\end{figure}

For the reconstruction scheme given by equation (\ref{eq:frink}), i.e., $\Delta U_{ij}=\alpha_d (W_{ij} - V_{ij})$, the semi-discrete scheme can be written as
\begin{equation}
|C_i| \dd{U_i}{t} = -\sum_{j \in N(i)} \left\{ P_{ij}  \frac{\theta_{ij} \Delta U_{ij}}{V_{ij} - U_i} (V_{ij} - U_i) + Q_{ij} \left[1 + \frac{\theta_{ji}  \Delta U_{ji}}{U_j - U_i} \right] (U_j - U_i) \right\}
\end{equation}
Defining
\begin{equation*}
r_{ij} = \frac{\Delta U_{ij}}{\frac{1}{2}(U_j - U_i)} = \frac{\alpha_d (W_{ij} - V_{ij})}{\frac{1}{2}(U_j - U_i)}
\end{equation*}
The above scheme can be written as
\begin{equation}
|C_i| \dd{U_i}{t} = -\sum_{j \in N(i)} \left\{ \frac{1}{2}P_{ij}  \theta_{ij} r_{ij} \frac{U_j - U_i}{V_{ij} - U_i} (V_{ij} - U_i) + Q_{ij} \left[1 - \frac{1}{2} \theta_{ji} r_{ji} \right] (U_j - U_i) \right\}
\end{equation}
The first term contains the reconstruction in cell $C_i$ while the second term is due to the reconstruction in cell $C_j$. From the coefficients of the terms containing $P_{ij}$ and $Q_{ji}$, this scheme will have positive coefficients provided the following two conditions are satisfied by the limiter function $\theta_{ij}$,
\begin{displaymath}
\theta_{ij} r_{ij} \frac{U_j - U_i}{V_{ij} - U_i} \le 0, \qquad 1 - \frac{1}{2} \theta_{ij} r_{ij} \ge 0
\end{displaymath}
The above two conditions will be satisfied if we choose the limiter as follows
\begin{equation}
\theta_{ij} = \begin{cases}
0 & (U_j - U_i)(U_i - V_{ij}) \le 0 \\
\theta(r_{ij}) & \textrm{otherwise}
\end{cases}
\label{eq:f1}
\end{equation}
where the function $\theta(r)$ must satisfy the following conditions:
\begin{enumerate}
\item $\theta(r)=0$ for all $r < 0$
\item $0 \le \theta(r) \le 1$ for all $r \in [0,\infty)$
\item $\theta(r) \le \frac{2}{r}$ for all $r > 0$
\end{enumerate}
These conditions are satisfied if we choose
\begin{equation}
\theta(r) = \max(0, \min(1, 2/r))
\label{eq:theta}
\end{equation}
The function $\theta(r)$ is sketched in figure~(\ref{fig:limiter}).  The first condition in equation (\ref{eq:f1}) checks whether $U_i$ lies between the values $V_{ij}$ and $U_j$; if it does not lie in this range, the limiter is set to zero and the order of the reconstruction is reduced to first order.
\subsection{Upwind scheme}

For the reconstruction scheme given by equation (\ref{eq:new}), i.e., $\Delta U_{ij}=\beta_d (U_i - V_{ij})$, the semi-discrete scheme can be written as
\begin{equation}
|C_i| \dd{U_i}{t} = \sum_{j \in N(i)} \left\{ P_{ij}  \theta_{ij} \beta_d (V_{ij} - U_i) - Q_{ij} \left[1 + \frac{\theta_{ji}  \Delta U_{ji}}{U_j - U_i} \right] (U_j - U_i) \right\}
\end{equation}
In this case, the first term inside the curly braces already satisfies the positivity condition. The coefficient in the second term is positive if we choose $\theta_{ij}=\theta(r_{ij})$, with $\theta(r)$ being given by equation~(\ref{eq:theta}) and
\begin{equation*}
r_{ij} = \frac{\Delta U_{ij}}{\frac{1}{2}(U_j - U_i)} = \frac{\beta_d (W_{ij} - V_{ij})}{\frac{1}{2}(U_j - U_i)}
\end{equation*}
We note in this case that there are less conditions on the limiter function $\theta_{ij}$ as compared to the previous scheme since we do not have to enforce the first condition in equation~(\ref{eq:f1}). This is due to the fact that the reconstruction scheme does not make use of face average values $W_{ij}$ but only the vertex and cell values.

\section{Numerical results for Euler equations}

We apply the schemes developed here to inviscid compressible flows governed by Euler equations of gas dynamics. These are a system of hyperbolic PDE representing the conservation of mass, momentum and energy. In three space dimensions, the vector of unknowns $U$ and the Cartesian components of the flux vector are given by
\begin{equation}
U=\begin{bmatrix}
\rho \\
 \rho u_1 \\
 \rho u_2 \\
 \rho u_3 \\
 E
\end{bmatrix},\quad
F_i = \begin{bmatrix}
\rho u_i \\
p\delta_{i1}+\rho u_1 u_i \\
p\delta_{i2}+\rho u_2 u_i \\
p\delta_{i3} + \rho u_3 u_i \\
 (E+p)u_i
\end{bmatrix}, \quad i=1, \ldots, 3
\end{equation}
where $\rho$ is the density, $(u_1,u_2,u_3)$ are the Cartesian components of the velocity, $p$ is the pressure and $E$ is the total energy per unit volume given by
\begin{equation}
E = \rho \varepsilon + \frac{1}{2} \rho |u|^2
\end{equation}
which is the sum of internal and kinetic energy. For an ideal gas, the pressure is related to the density and internal energy per unit mass $\varepsilon$  by
\begin{equation}
p = (\gamma-1) \rho \varepsilon
\end{equation}
where $\gamma$ is the ratio of specific heats and is a constant for a given gas. In all the numerical examples, we take $\gamma=1.4$ which corresponds to air under normal conditions.

The Euler equations are solved on tetrahedral grids using a matrix-free LUSGS scheme~\cite{sharov1998}. In the case of unsteady problems, the 3-stage strong stability preserving Runge-Kutta scheme is used~\cite{Shu1988439}. As in \cite{frink1994}, we reconstruct the primitive variables $(\rho, u_1, u_2, u_3, p)$ rather than the conserved variables. This is beneficial in maintaining the positivity of density and pressure. The numerical flux function is either based on the Roe scheme~\cite{Roe1981357} or the KFVS scheme~\cite{Mandal1994447}. The latter scheme is known to be entropy consistent~\cite{springerlink:10.1007/PL00001540} and hence very robust; we use it for the high Mach number problems which contain strong shocks. Most of the grids used in this work were generated using the open source tool GMSH~\cite{NME:NME2579}.

We also compare the current schemes with the limited schemes of Jameson~\cite{jameson-iccfd2000} where a central limiter is used. The reconstructed values on any face $S_{ij}$ are given by
\begin{equation}
U_{ij}^+ = U_i + \frac{1}{3} L(U_i - V_{ij}, V_{ji}-U_j), \quad
U_{ij}^- = U_j - \frac{1}{3} L(U_i - V_{ij}, V_{ji}-U_j)
\label{eq:jam}
\end{equation}
where $L(a,b)$ is a limited average which is taken to be
\begin{equation}
L(a,b) = \frac{1}{2}(a+b)[1-R(a,b)], \qquad R(a,b) = \left| \frac{a-b}{\max(|a|+|b|,\epsilon h^{3/2})} \right|^q, \quad q \in [1,3]
\label{eq:jam2}
\end{equation}
If $q=1$, then the above limiter reduces to the minmod limiter. Note that larger values of $q$ lead to less limiting. For this scheme, we are not able to prove the maximum stability property. Note that the term {\em Jameson scheme} is used to refer to the reconstruction scheme given by equations~(\ref{eq:jam}) and (\ref{eq:jam2}).

All the results are obtained using the modified Shepard interpolation. For some problems, a few of the vertices have negative weights which we do not modify. In each case, we indicate the number of negative weights using the original interpolation scheme of Frink and the modified scheme.
\subsection{Shock tube problem}

The standard shock tube problem is solved on a three dimensional grid. The computational domain is in the form of a channel with a square cross-section and whose axis is along the $x$-axis. The side walls of the channel are treated as slip walls. The computational grid contains about 100 points along the axis of the channel. The solution is taken along the center line of the channel for comparison with the exact solutions of the Riemann problem. Two standard Riemann problems are considered, the Sod problem and a problem which involves low densities in the solution. In both cases, the initial discontinuity is placed at $x=0.5$. The initial conditions defining the Riemann problem are shown in table~(\ref{tab:stube}). For Test 1 we use the Roe flux function while for Test 2, we use the KFVS flux.
\begin{table}
\begin{center}
\begin{tabular}{|l||c|c|c||c|c|c|}
\hline
Case & $\rho_l$ & $u_l$ & $p_l$ & $\rho_r$ & $u_r$ & $p_r$ \\
\hline
Test 1 & 1 & 0 & 1 & 0.125 & 0 & 0.1 \\
\hline
Test 2 & 1 & -2 & 0.4 & 1 & 2 & 0.4 \\
\hline
\end{tabular}
\end{center}
\caption{Initial conditions for shock tube problem}
\label{tab:stube}
\end{table}
Limited versions of the reconstruction schemes are used since the solutions have shock and contact discontinuities. With the original interpolation scheme, 91 vertices have some negative weights with the smallest weight being -0.03, while the smallest determinant is $1.44 \times 10^{-15}$. With the modified interpolation scheme, there is only one vertex with a negative weight of $-1.5 \times 10^{-5}$ while the smallest determinant is 1.18.

\begin{figure}
\begin{center}
\includegraphics[width=0.50\textwidth]{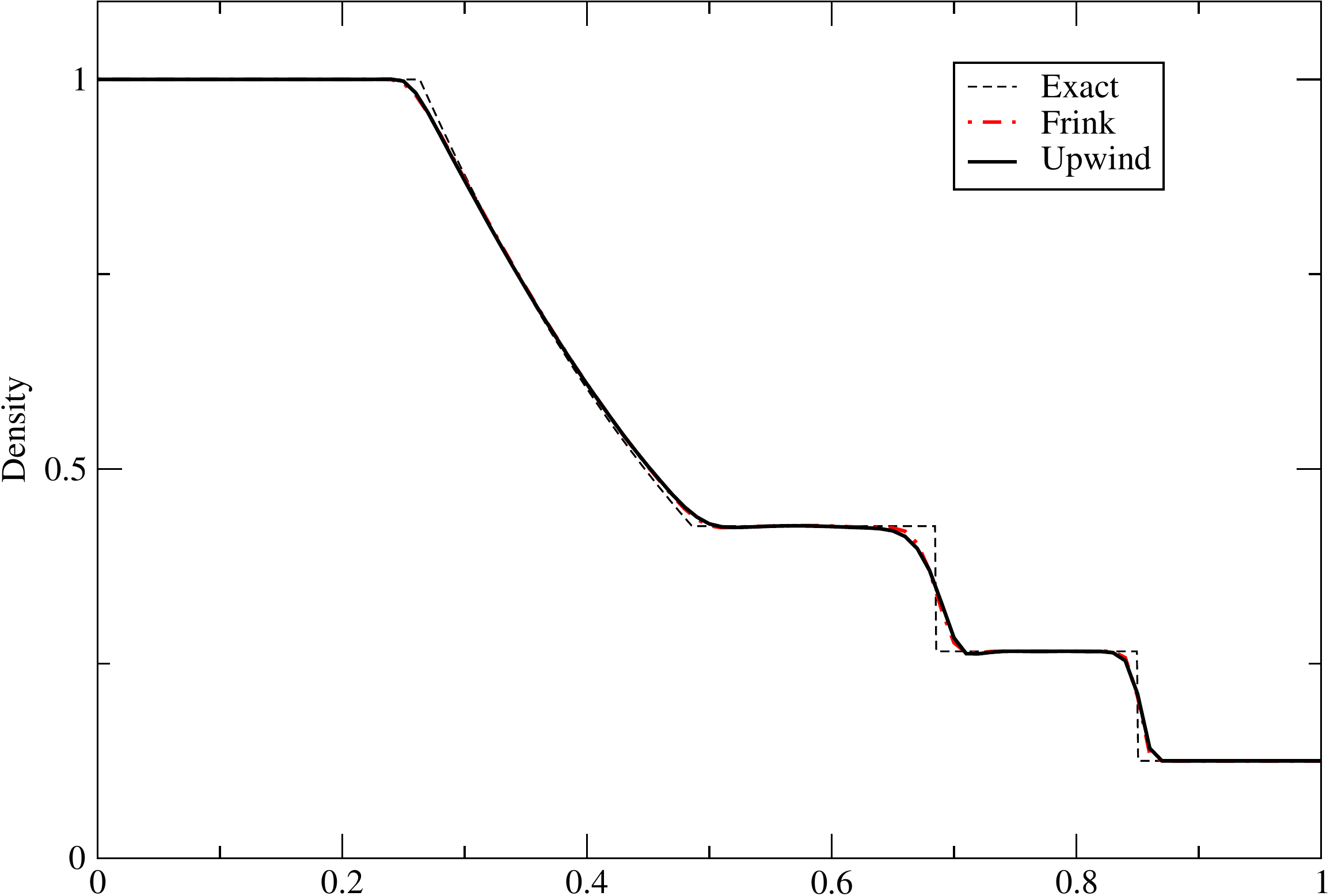}
\includegraphics[width=0.35\textwidth]{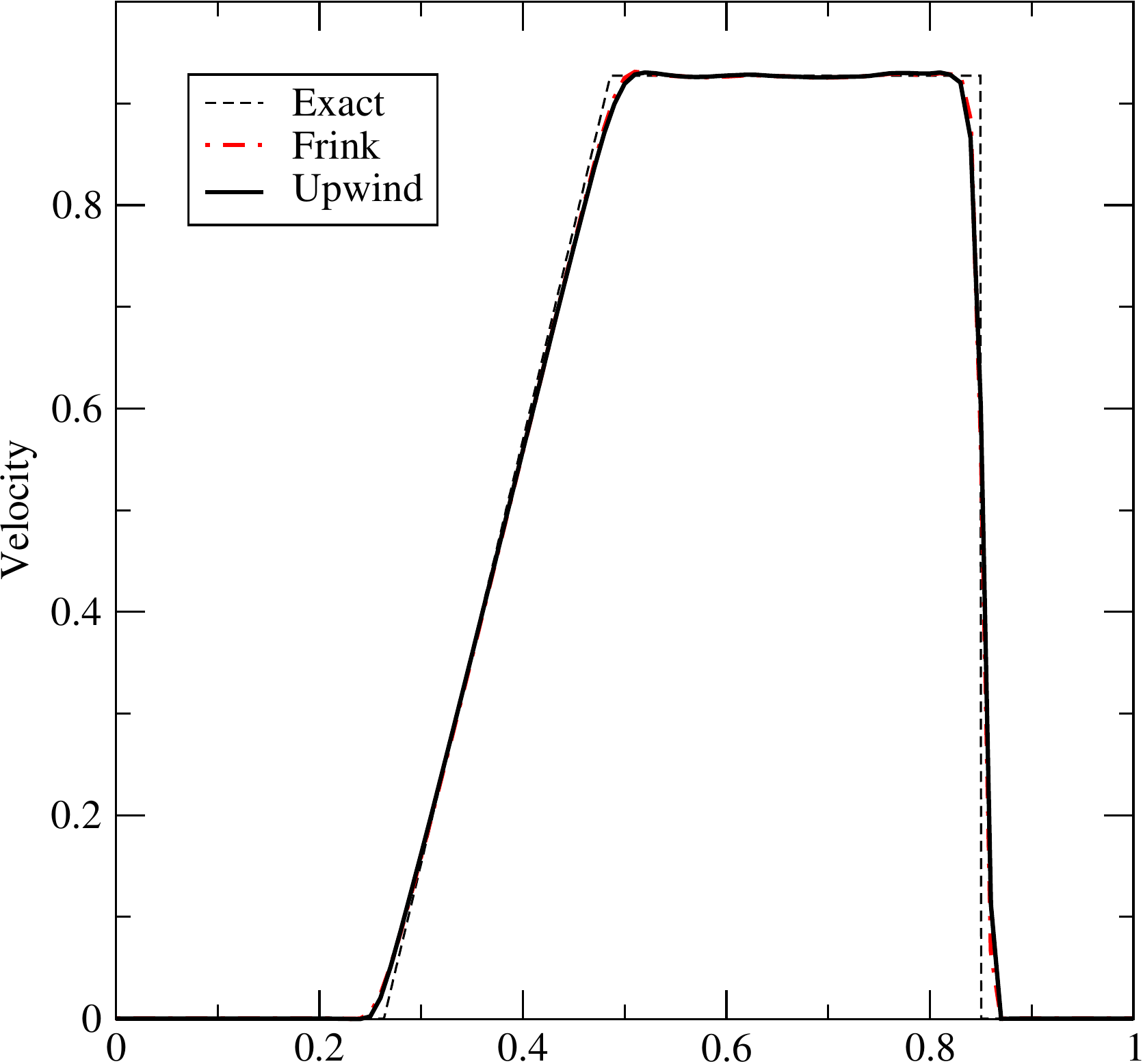}
\end{center}
\caption{Solution for Test case 1 at time $t=0.2$}
\label{fig:test1}
\end{figure}

The scheme given by equation~(\ref{eq:jam}) did not work for Test 1. Figure~(\ref{fig:test1}) shows the density and velocity obtained for Test 1; all the other schemes are able to give oscillation free solutions. There is no appreciable difference between the Frink and upwind reconstructions. 

\begin{figure}
\begin{center}
\begin{tabular}{ccc}
\includegraphics[width=0.32\textwidth]{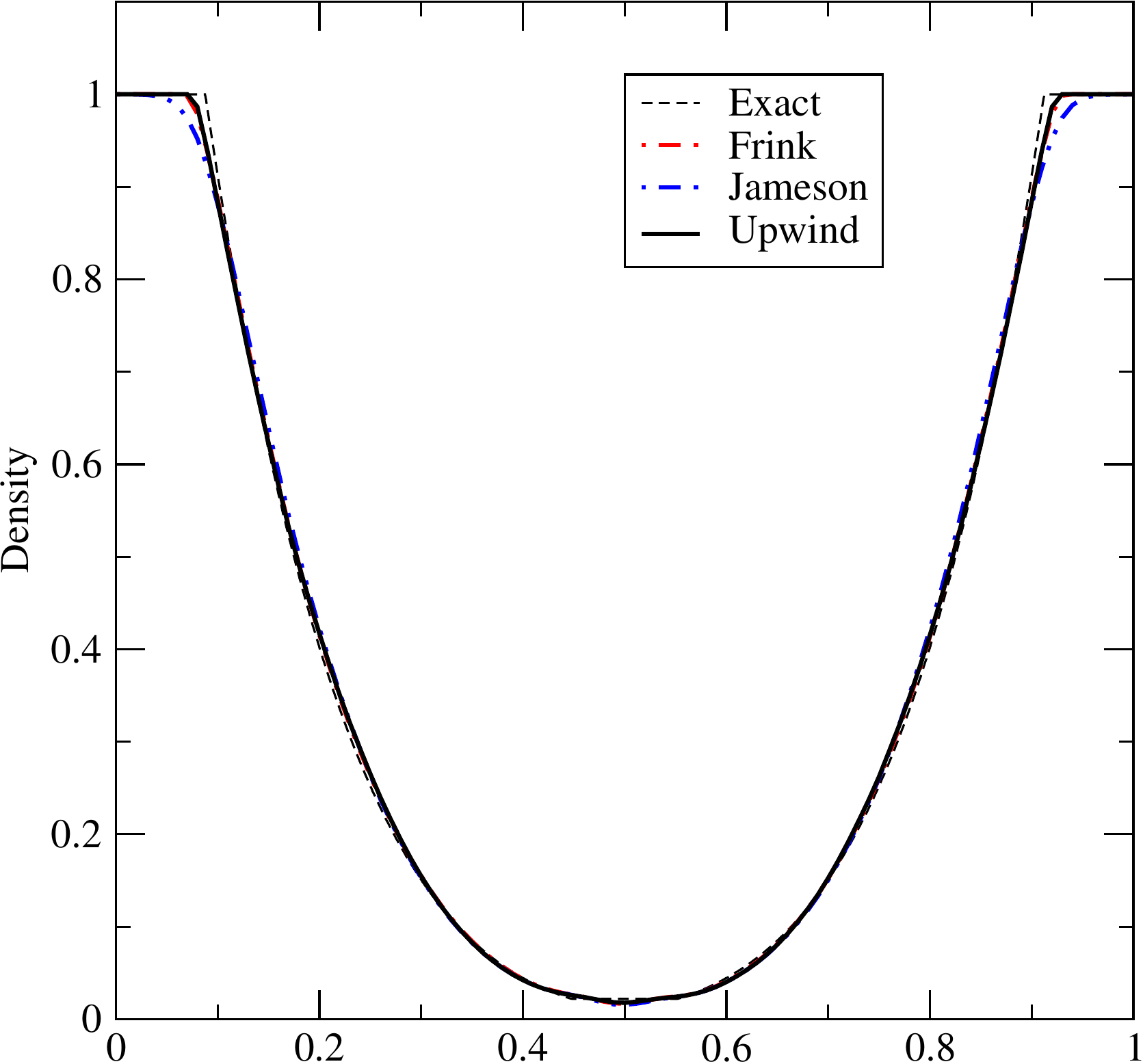} &
\includegraphics[width=0.32\textwidth]{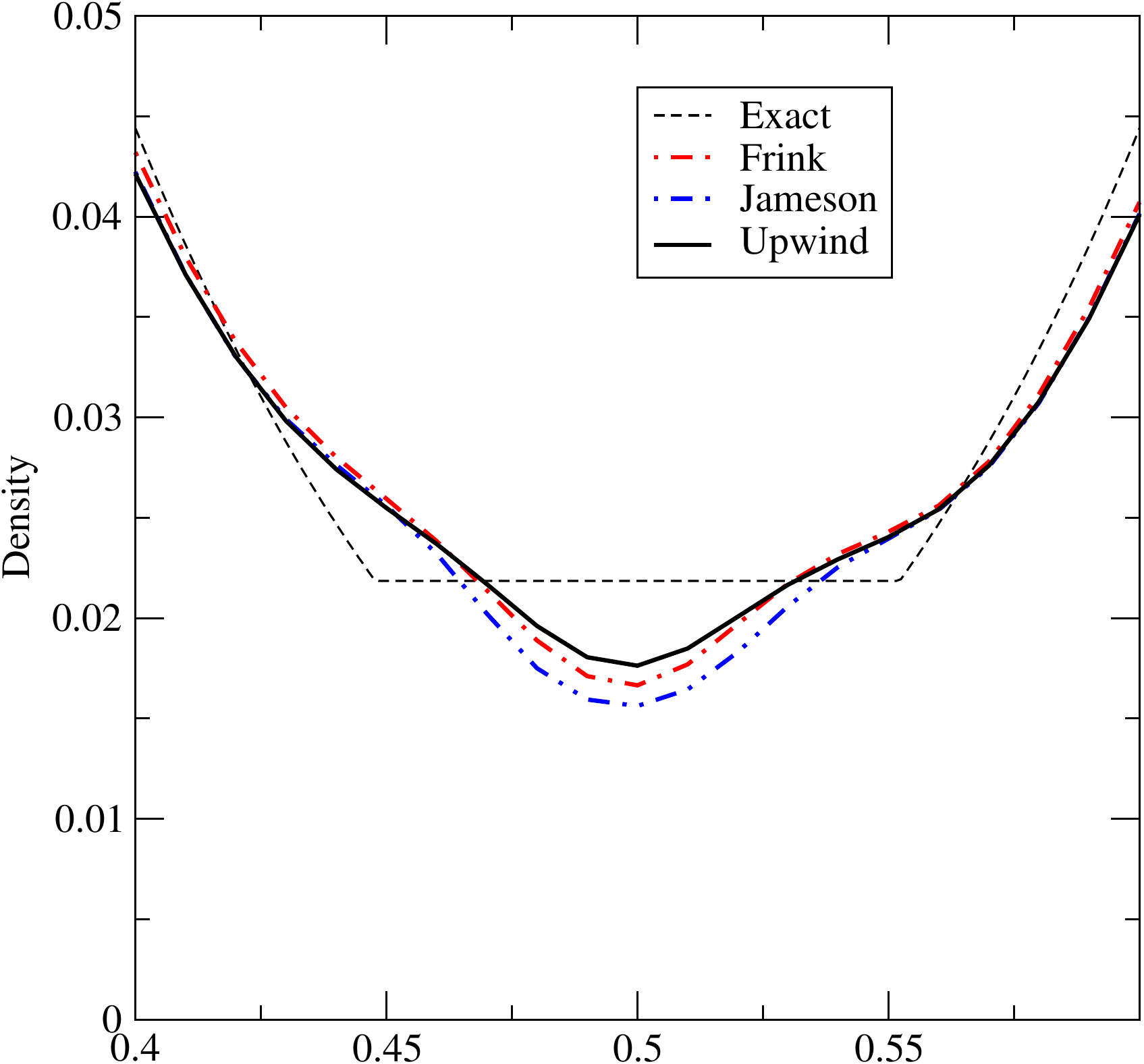} &
\includegraphics[width=0.32\textwidth]{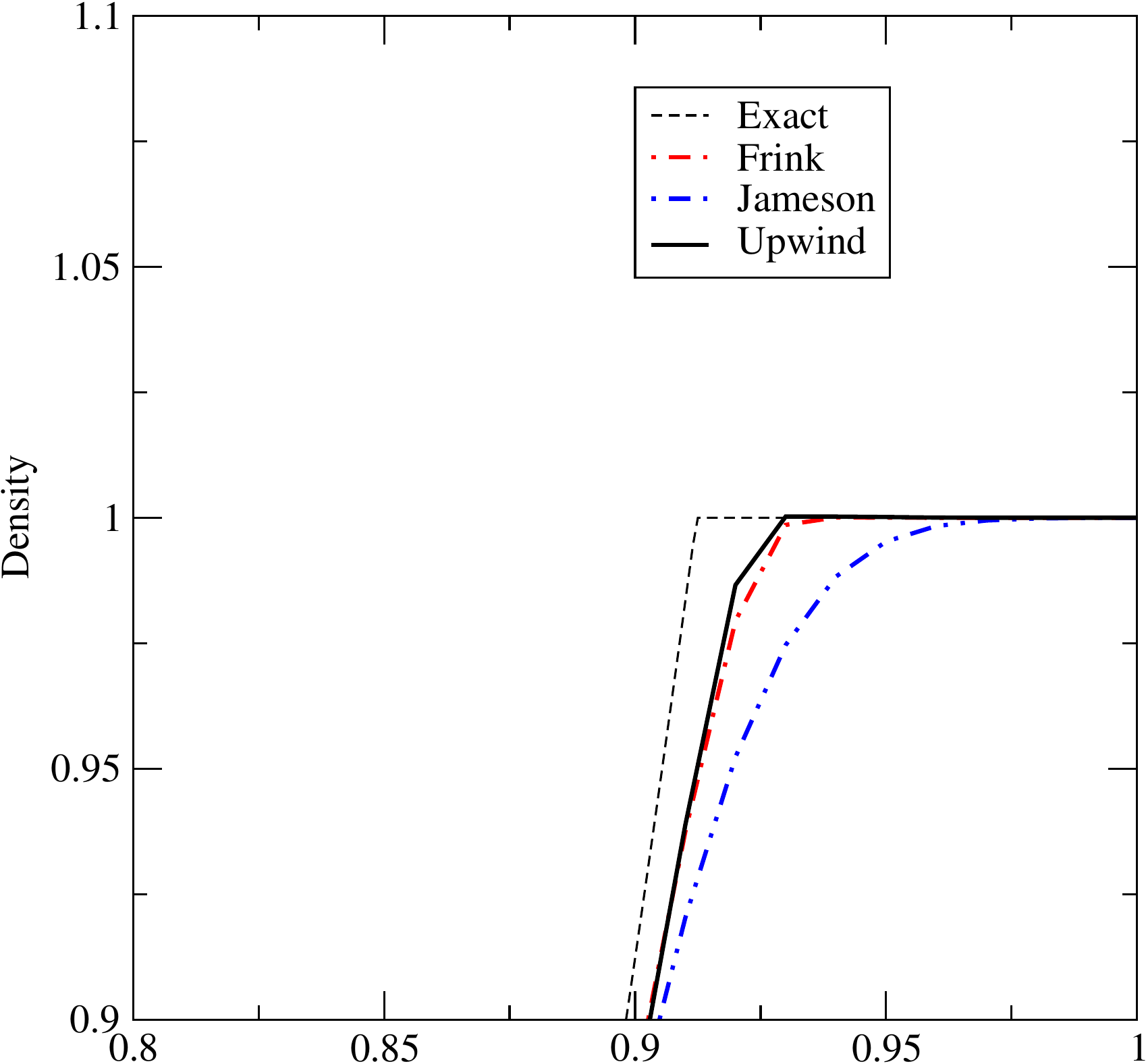} \\
(a) & (b) & (c) \\
\end{tabular}
\end{center}
\caption{Solution for Test case 2 at time $t=0.15$}
\label{fig:test2}
\end{figure}

Figure~(\ref{fig:test2}) shows the density for Test 2 obtained from all the three schemes. Due to the expansion wave, the density reaches a very low value at the middle of the computational domain as shown in the zoomed view in figure~(\ref{fig:test2}-b). All the schemes are able to preserve the positivity of density. The upwind scheme has the least overshoot of density while Jameson scheme has the highest. The zoomed view in figure~(\ref{fig:test2}-c) show the tip of the expansion wave where the solution has a corner shape. We see that the two new limited reconstructions proposed here give the best results in terms of preserving the corner shape.

\subsection{Transonic flow over Onera M6 wing}
This is a very standard test case for compressible flows using finite volume methods for which experimental results are available. The problem involves transonic flow over the Onera M6 wing at a free stream Mach number of 0.839 and an angle of attack of 3.06 degrees~\cite{agard-ar-138}. The solution consists of shocks on the upper surface of the wing with a lambda structure on the wing. The tetrahedral grid used in the computations consists of 341797 cells. With the original vertex averaging procedure, there are 2044 vertices with negative weights and the minimum determinant is 1.15E-15, while with the modified averaging procedure, there are only 86 vertices with negative weights and the minimum determinant is 0.098.

\begin{figure}
\begin{center}
\begin{tabular}{cc}
\includegraphics[width=0.48\textwidth]{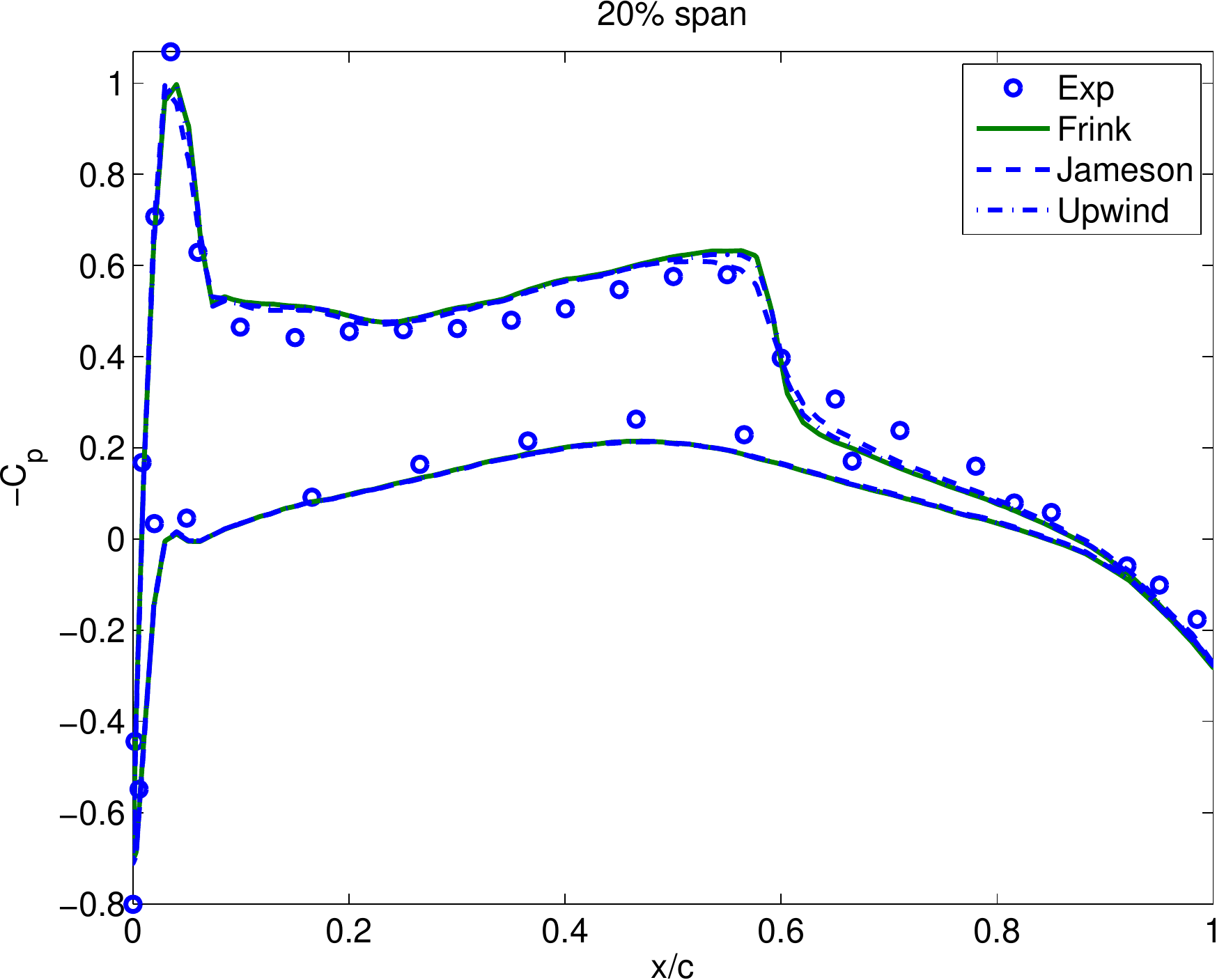} &
\includegraphics[width=0.48\textwidth]{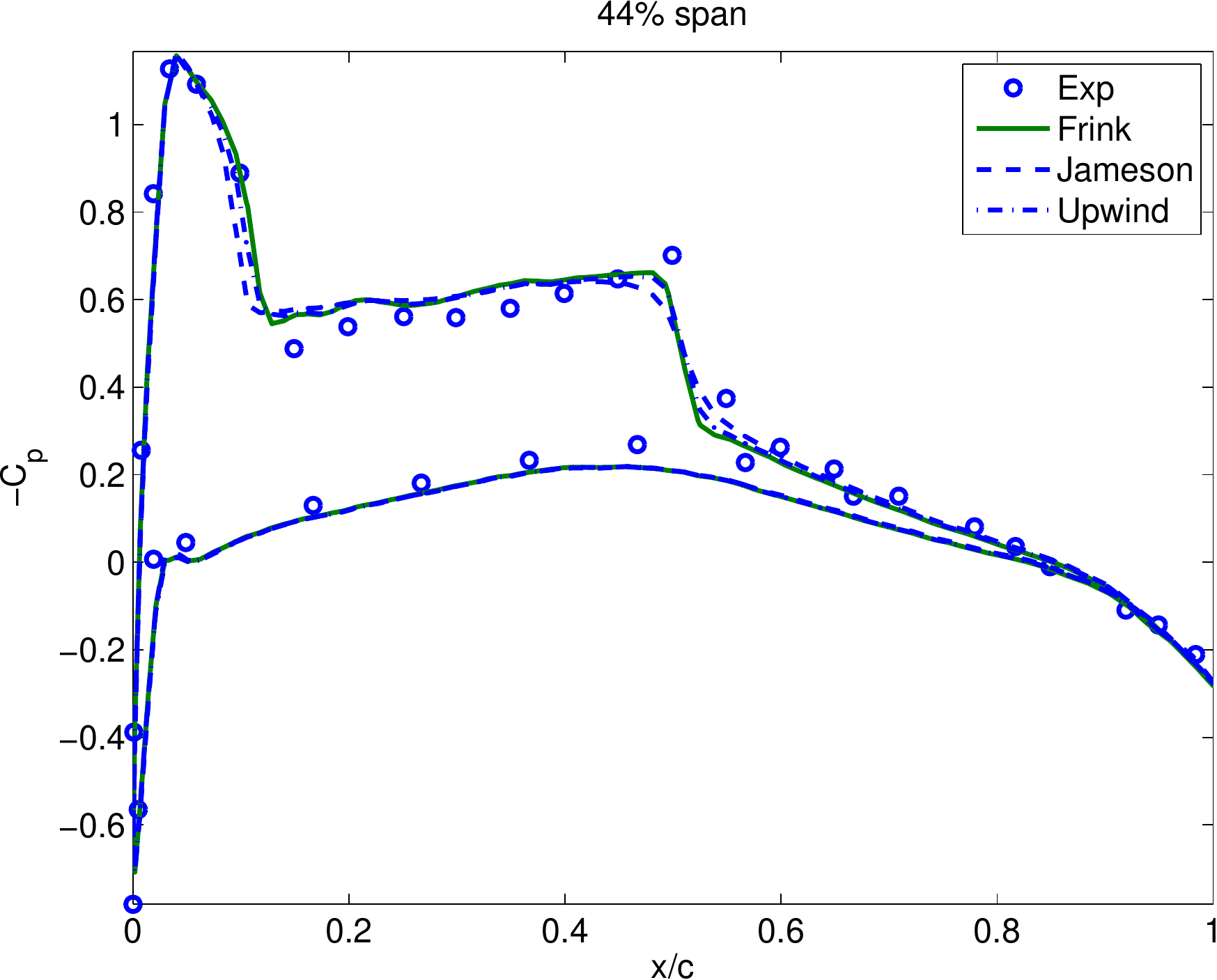} \\
\includegraphics[width=0.48\textwidth]{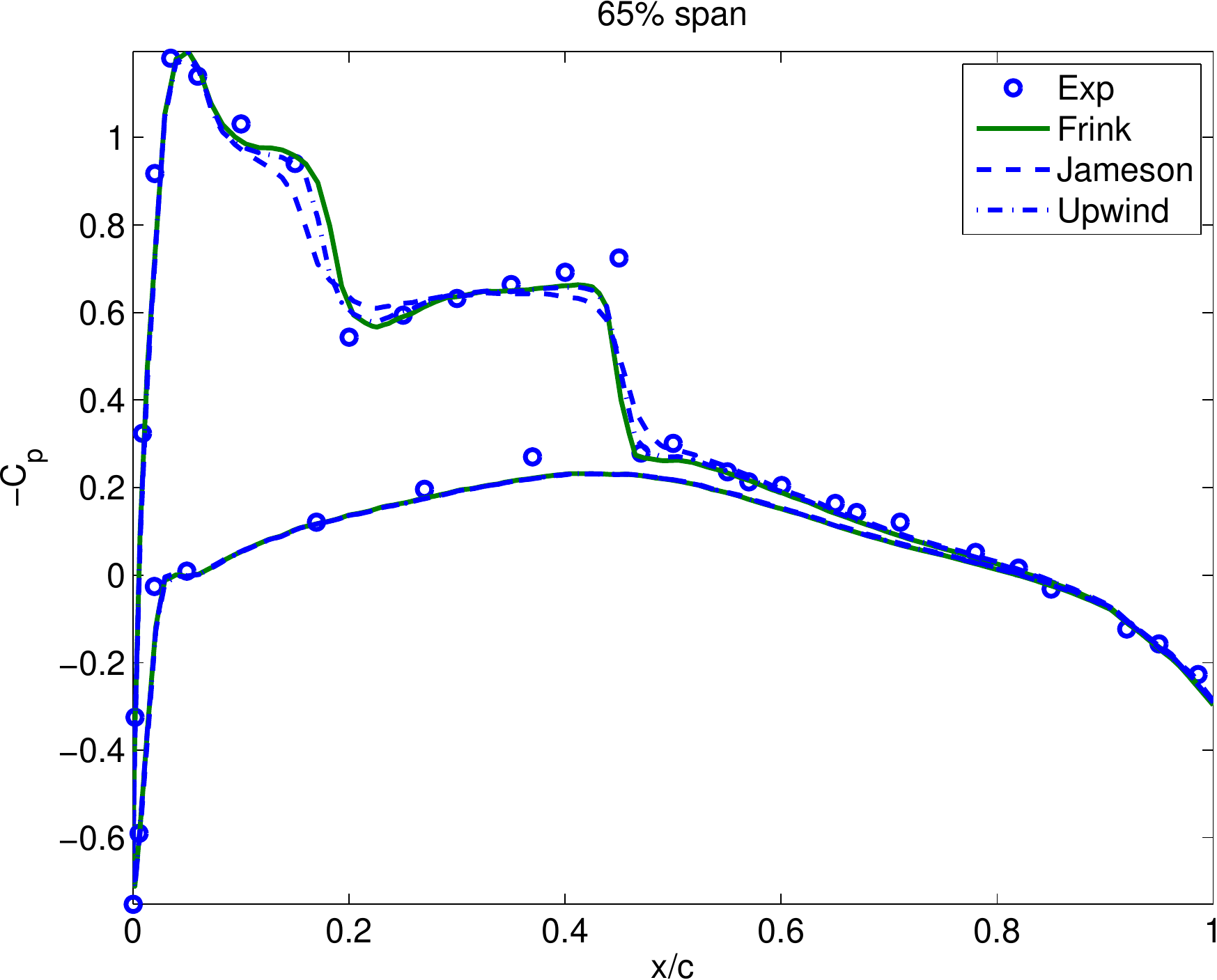} &
\includegraphics[width=0.48\textwidth]{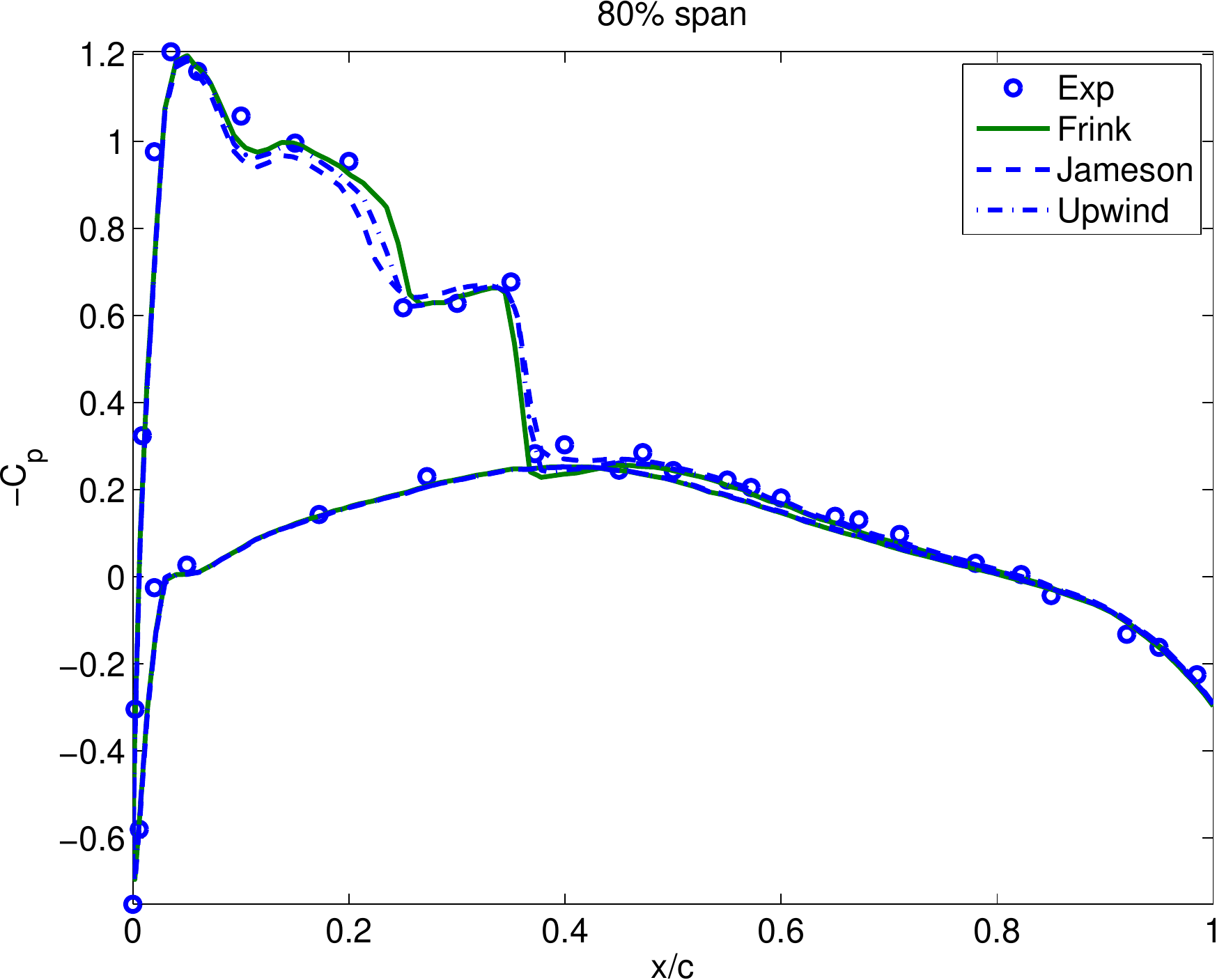} \\
\end{tabular}
\end{center}
\caption{Pressure coefficient for Onera M6 wing}
\label{fig:oneracp}
\end{figure}

\begin{figure}
\begin{center}
\begin{tabular}{ccc}
\includegraphics[width=0.33\textwidth]{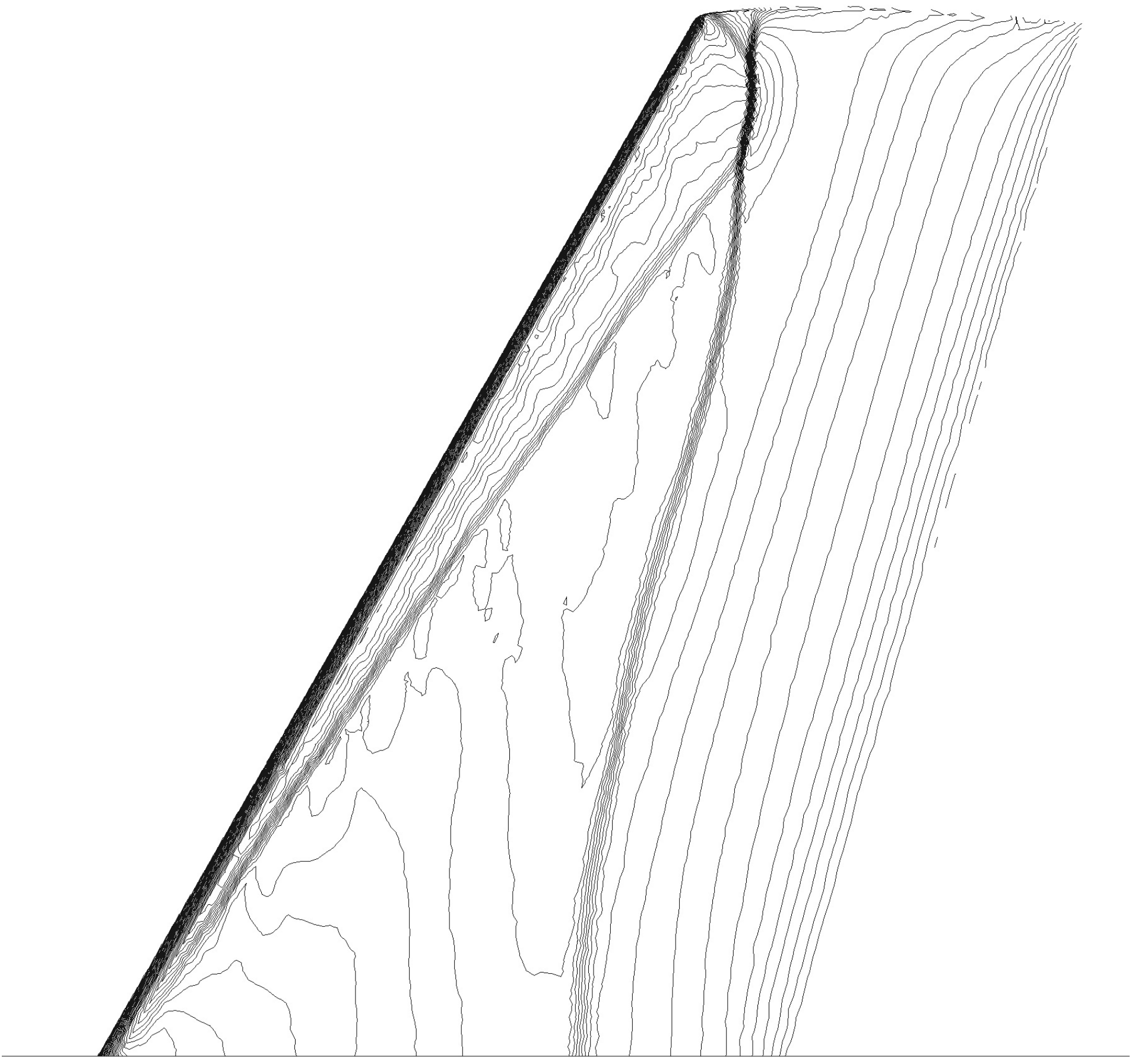} &
\includegraphics[width=0.33\textwidth]{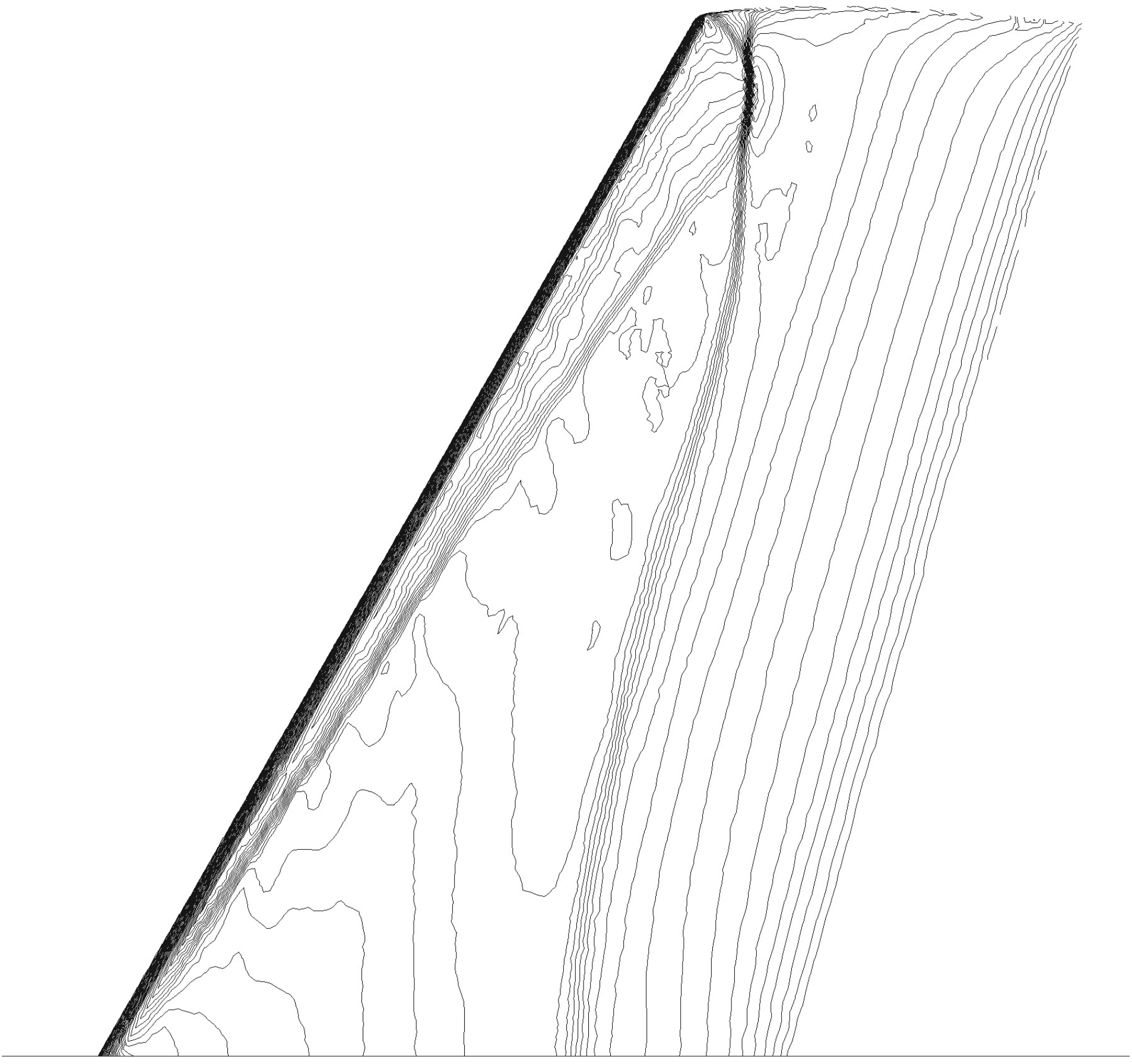} &
\includegraphics[width=0.33\textwidth]{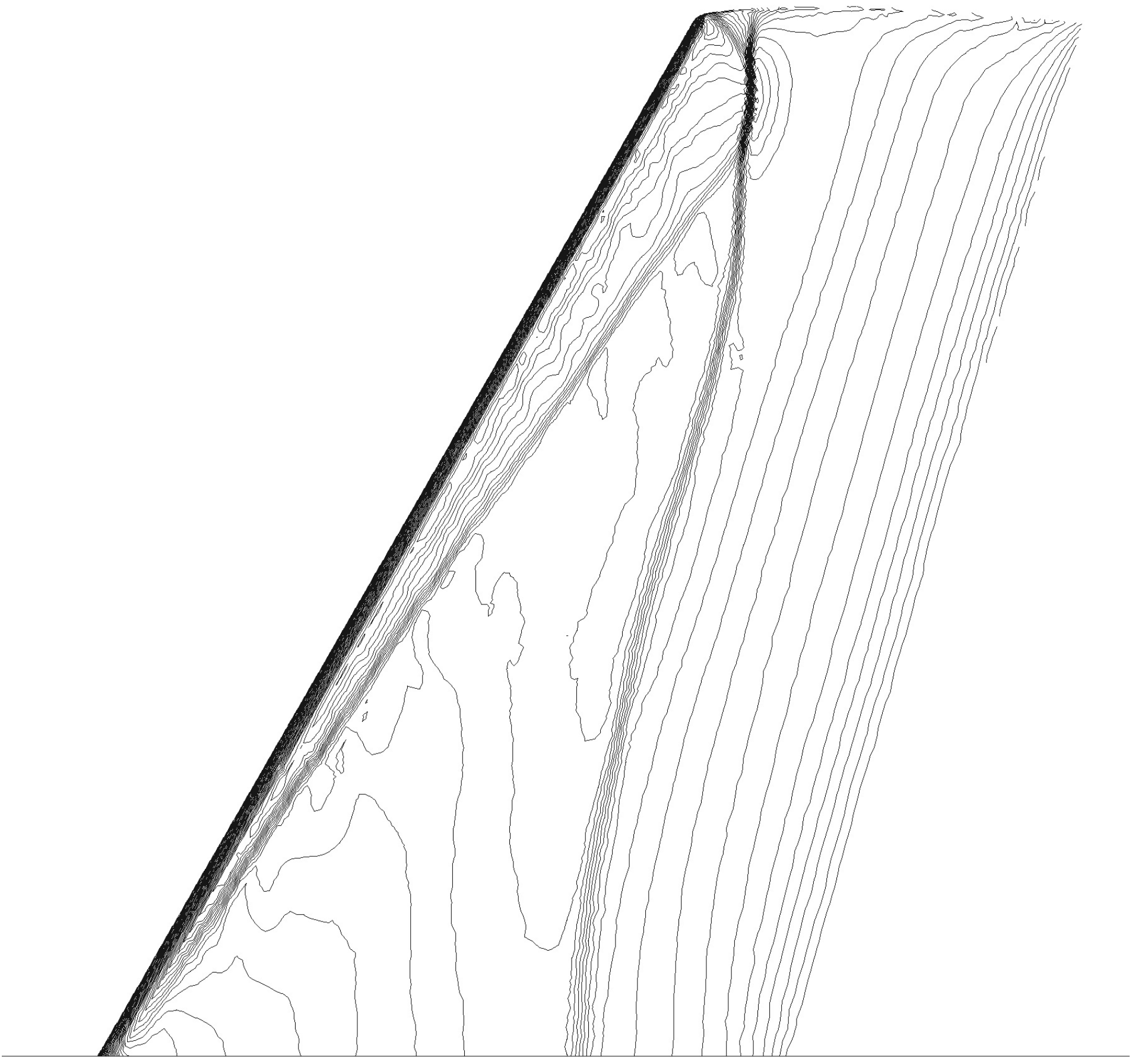} \\
Frink & Jameson & Upwind \\
\end{tabular}
\end{center}
\caption{Pressure contours for Onera M6 wing}
\label{fig:onerapre}
\end{figure}

The pressure coefficient at different spanwise stations are plotted and compared with experimental results in figure~(\ref{fig:oneracp}). The Frink and upwind reconstruction schemes yield very similar results, while the Jameson scheme is slightly dissipative at the shocks. The comparison with experimental results is remarkably good for all the schemes considering that the grid used is not very fine. The pressure variation on the wing surface is shown in figure~(\ref{fig:onerapre}) for the three schemes with the same contour levels; the lambda shock is well captured by all of them. However, the Jameson reconstruction can be again seen to give slightly more diffused shock compared to the other two schemes.

\subsection{Supersonic flow over hemi-sphere}

This problem consists of a Mach 3 flow over a hemi-sphere which leads to the formation of a bow shock in front of the sphere. We use the KFVS scheme for the numerical flux function, which is known to be entropy consistent and gives physically correct solutions even for strong shocks. The radius of the hemi-sphere is one and the problem is solved on two grids: grid G1 has a cell size of 0.05 on the surface of the hemi-sphere with 908330 tetrahedra while grid G2 has a cell size of 0.025 with 1934791 tetrahedra. All the averaging weights were positive using both methods. In the original procedure of Frink, the smallest determinants were 2.9e-11 and 3.5e-13. In the case of the modified averaging procedure, the smallest determinants for the two grids were 1.54 and 1.98 respectively.

\begin{figure}
\begin{center}
\begin{tabular}{cc}
\includegraphics[width=0.45\textwidth]{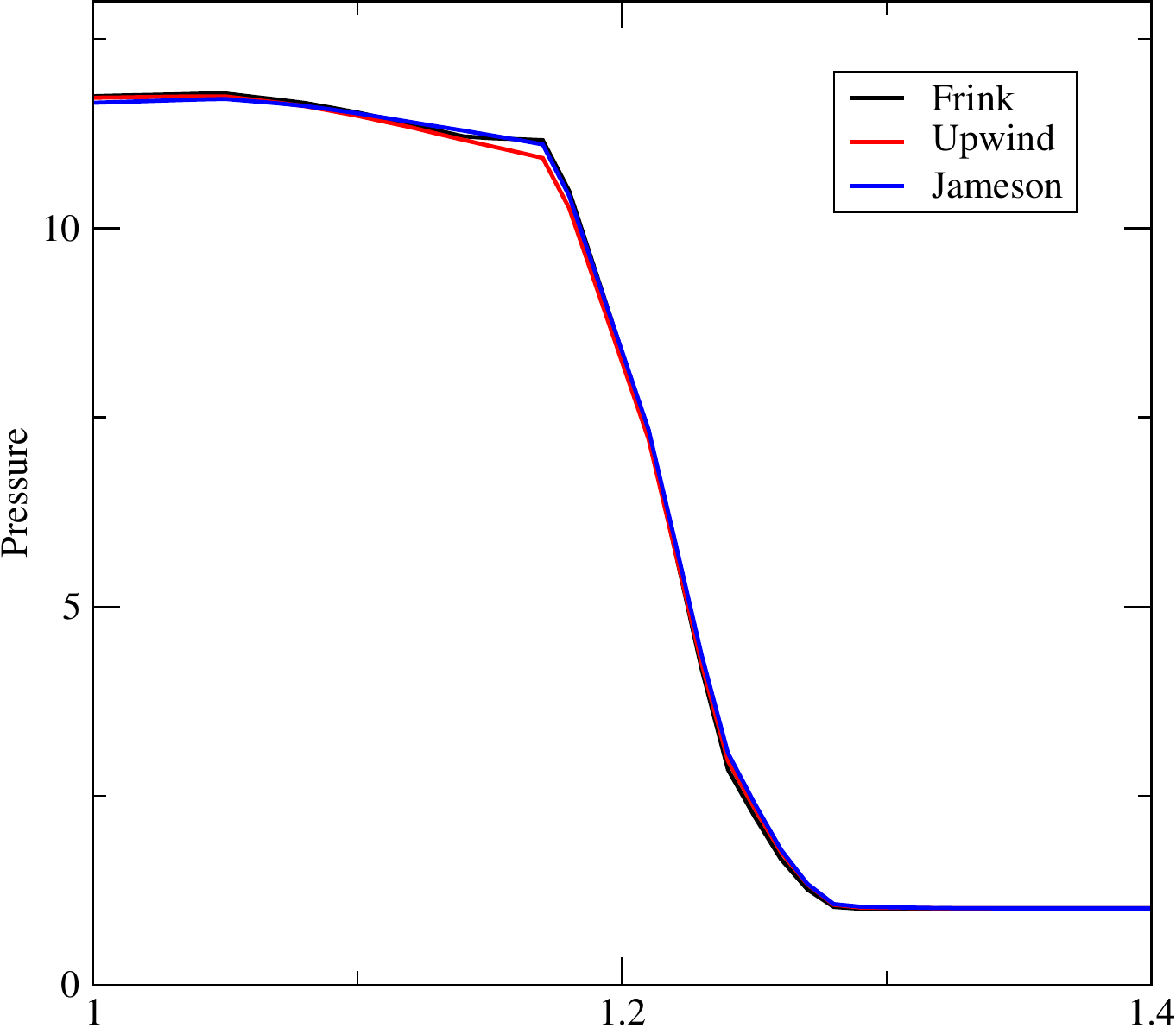} &
\includegraphics[width=0.45\textwidth]{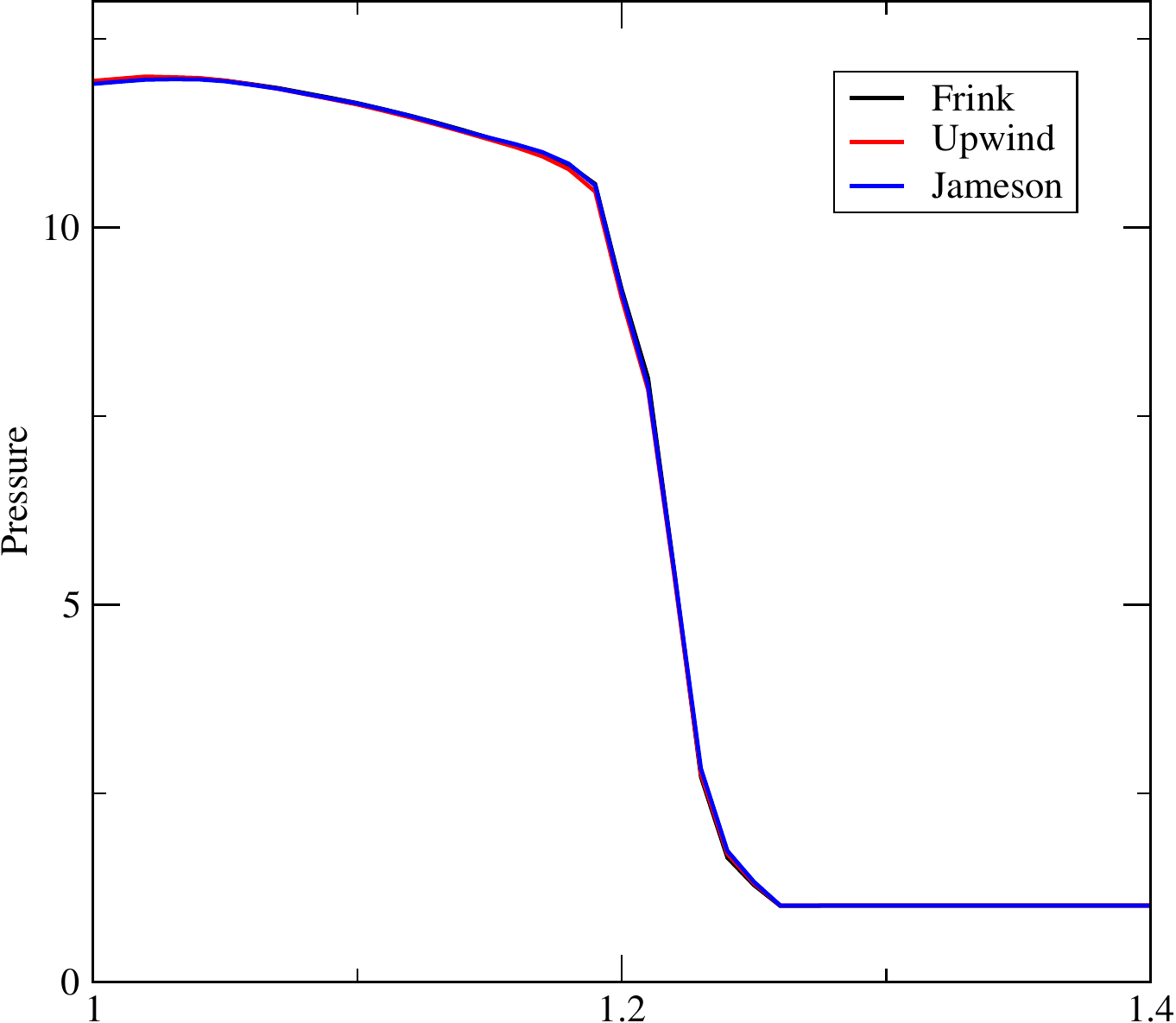} \\
(a) Grid G1 & (b) Grid G2
\end{tabular}
\caption{Pressure variation along stagnation stream-line for flow over hemi-sphere}
\label{fig:hemipline}
\end{center}
\end{figure}

\begin{figure}
\begin{center}
\begin{tabular}{ccc}
\includegraphics[width=0.4\textwidth,angle=-90]{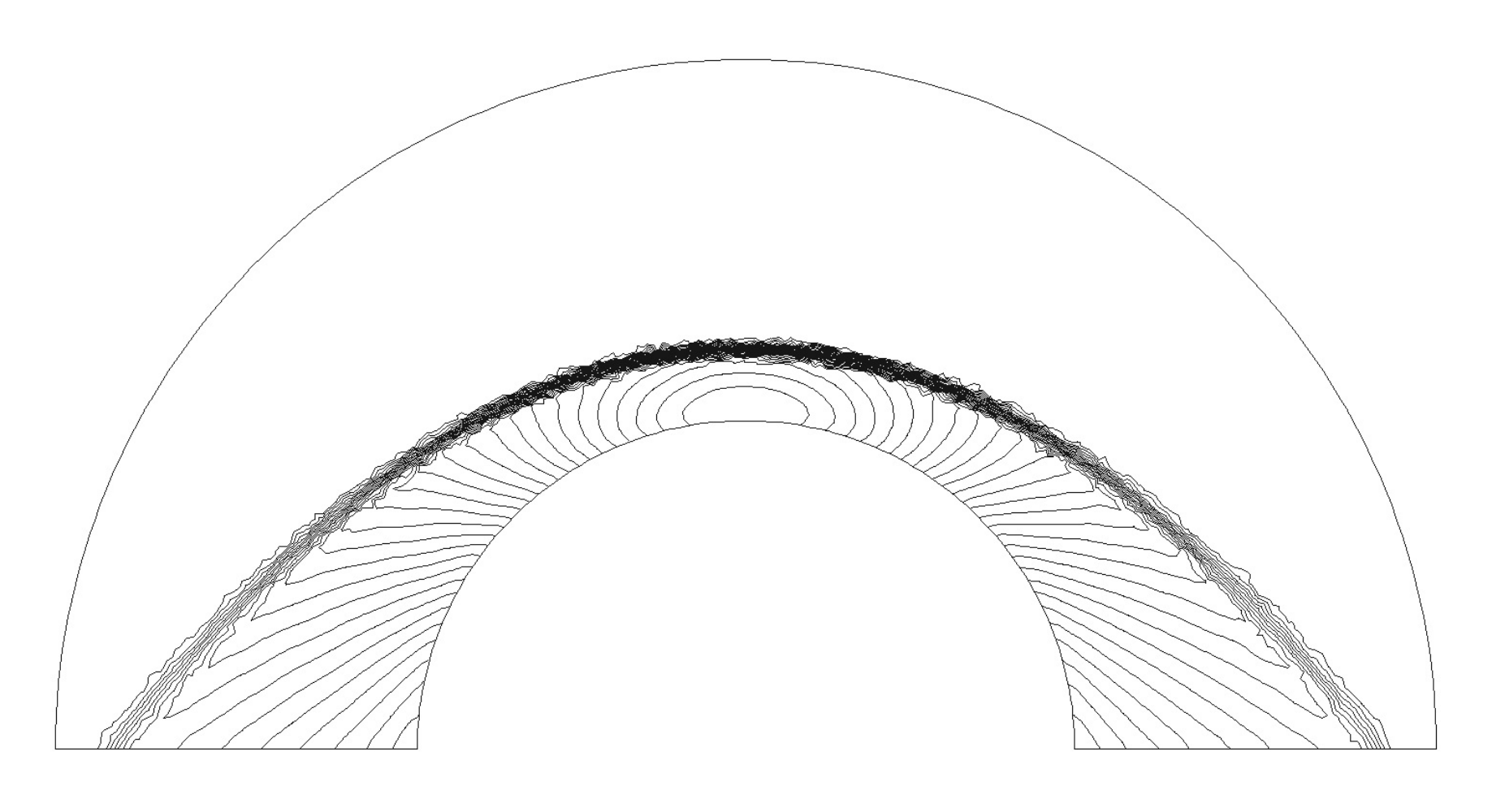} &
\includegraphics[width=0.4\textwidth,angle=-90]{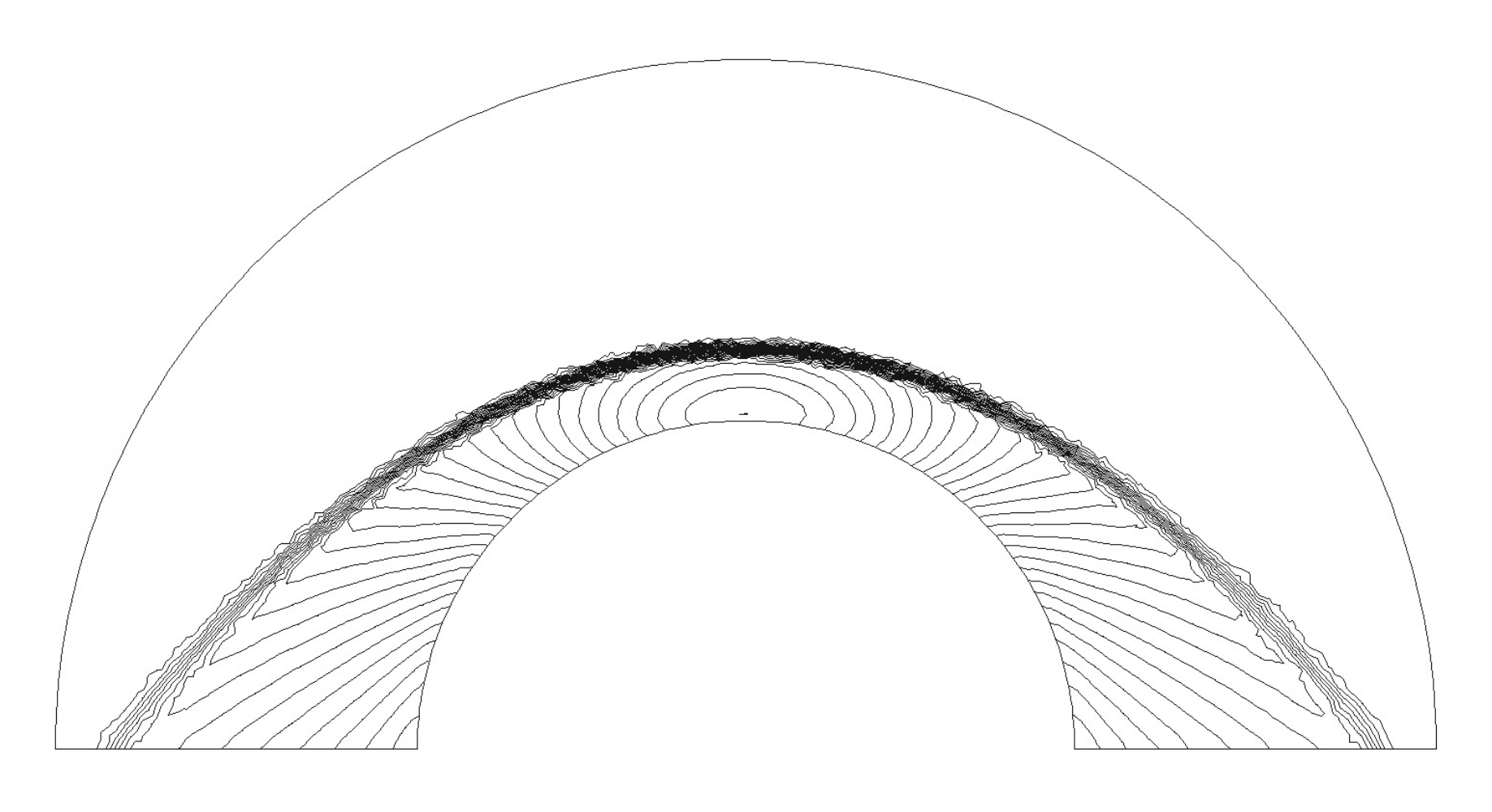} &
\includegraphics[width=0.4\textwidth,angle=-90]{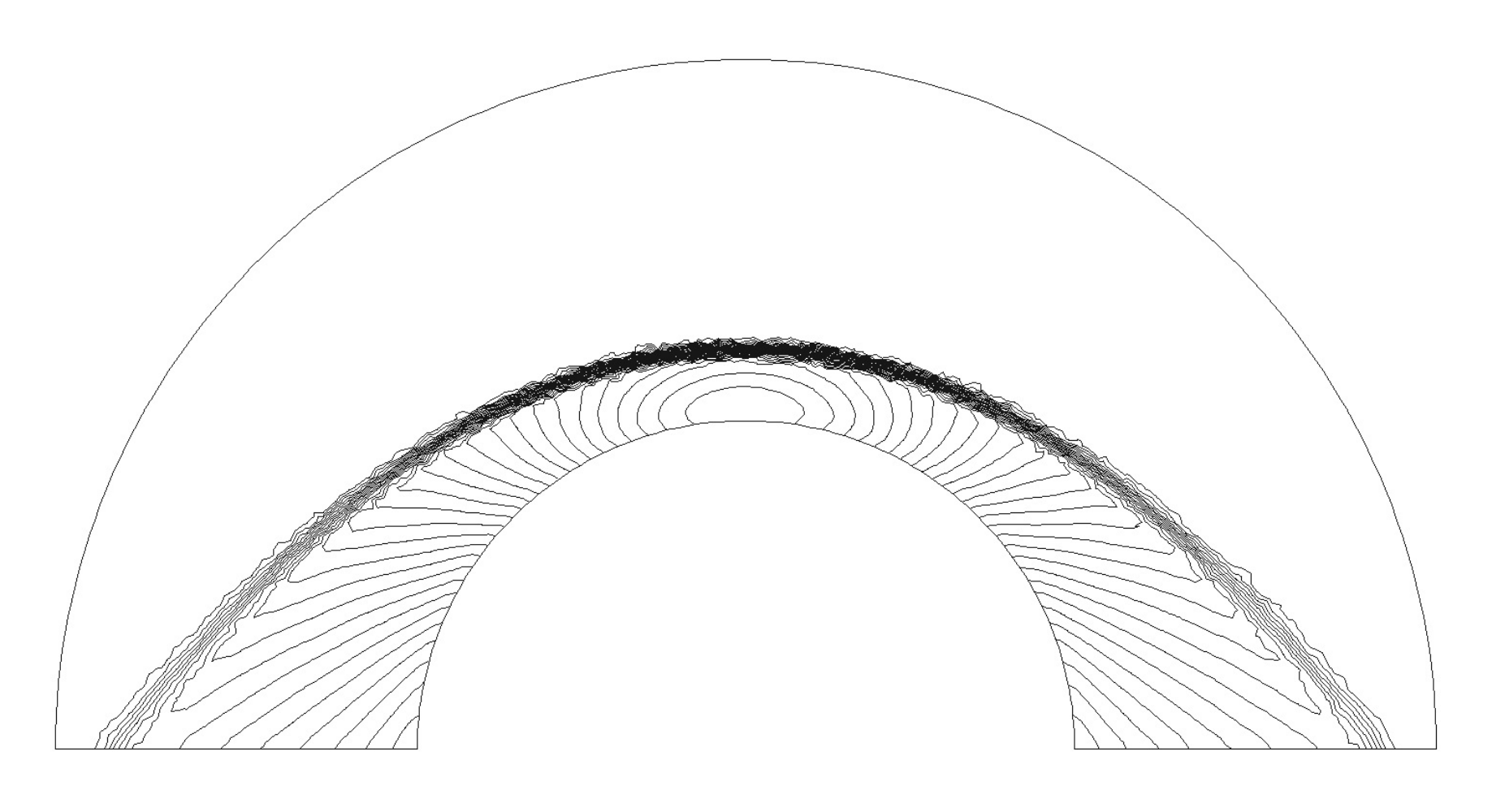} \\
(a) Frink & (b) Upwind & (c) Jameson
\end{tabular}
\caption{Pressure contour for flow over hemi-sphere using grid G2}
\label{fig:hemip}
\end{center}
\end{figure}

The variation of pressure along the stagnation stream-line is shown in figure~(\ref{fig:hemipline}). Along this line, the shock is normal to the line and we see that there are no oscillations in the pressure distribution. This is true of the density and velocity also but the plots are not shown here. The shock width is seen to decrease when we go from the coarser grid G1 to the finer grid G2. Figure~(\ref{fig:hemip}) shows the pressure contours obtained on G2 using the three reconstruction methods. It is clear that all of them give similar results without any oscillations.

\subsection{Blast wave}

In order to test the limiters under strong shocks, we consider a problem with a large temperature/pressure discontinuity. The initial condition consists of a small region of high temperature. The solution develops into an expanding shock wave. In the initial condition, the density is uniform with a value of 1.228 while the velocity is zero. The temperature inside a spherical core of radius 5 meters is $8.1 \times 10^7$ K while outside the core it is 298 K. The numerical flux function is based on KFVS scheme. The computational domain is a cube whose sides are 81 meters and the grid consists of 3182815 tetrahedra. With the original averaging procedure, there was one vertex with negative weight and the minimum determinant was $8.7 \times 10^{-4}$. With the modified averaging procedure, there was no negative weight and the minimum determinant was 1.5.

\begin{figure}
\begin{center}
\begin{tabular}{cc}
\includegraphics[width=0.48\textwidth]{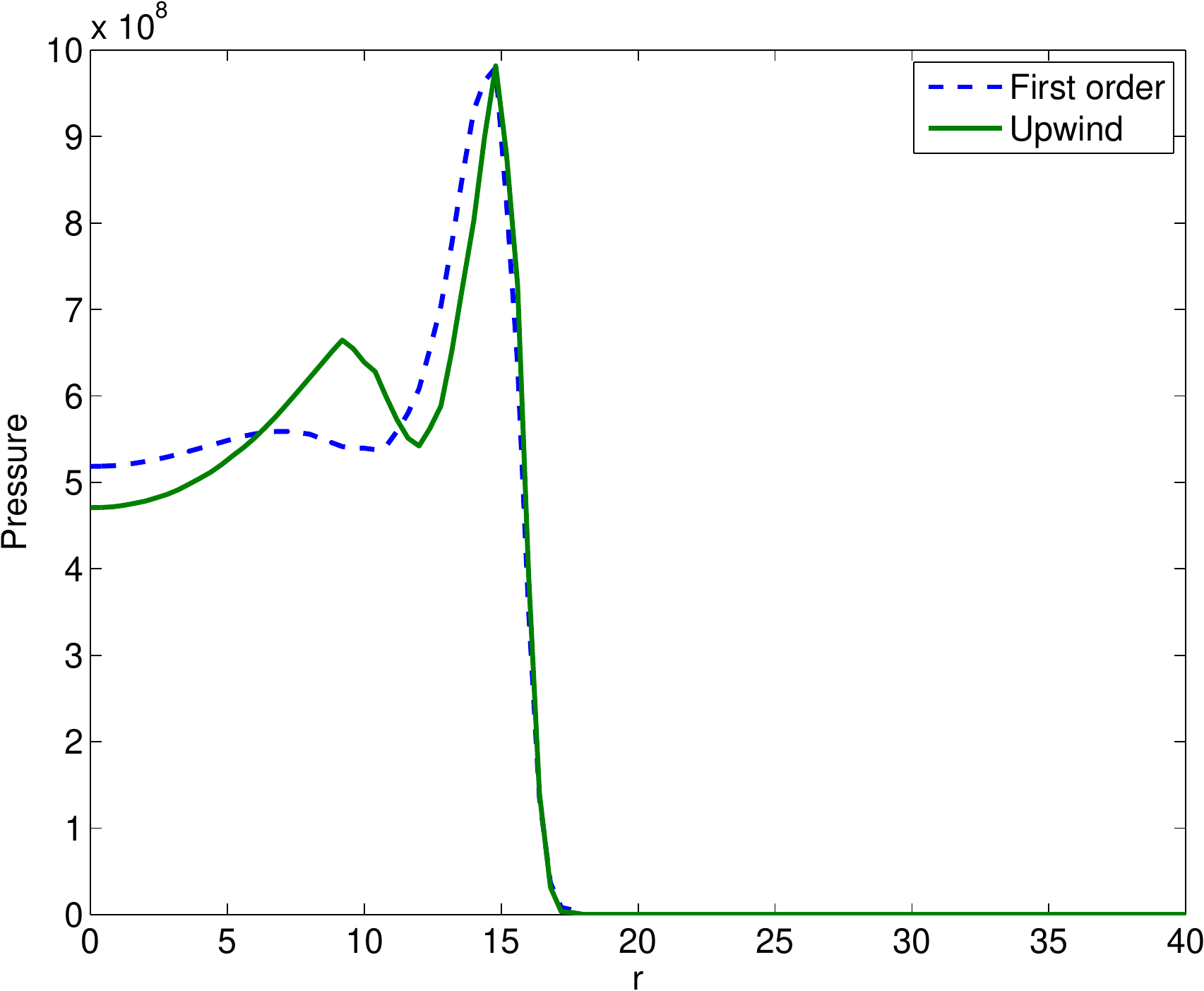} &
\includegraphics[width=0.48\textwidth]{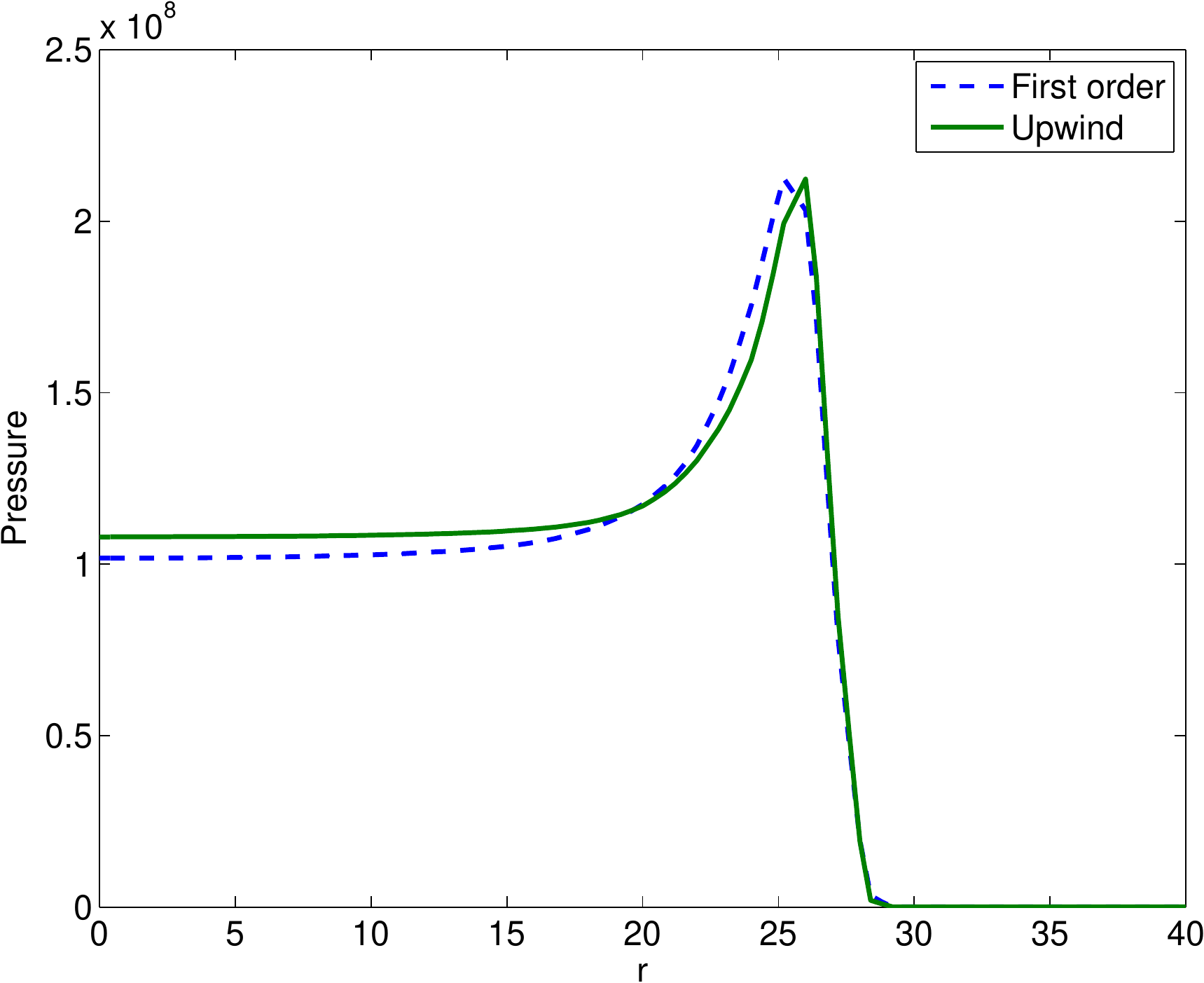} \\
(a) 25 time steps & (b) 100 time steps
\end{tabular}
\caption{Pressure variation along a radial line for blast wave problem}
\label{fig:blastp}
\end{center}
\end{figure}

\begin{figure}
\begin{center}
\includegraphics[width=0.48\textwidth]{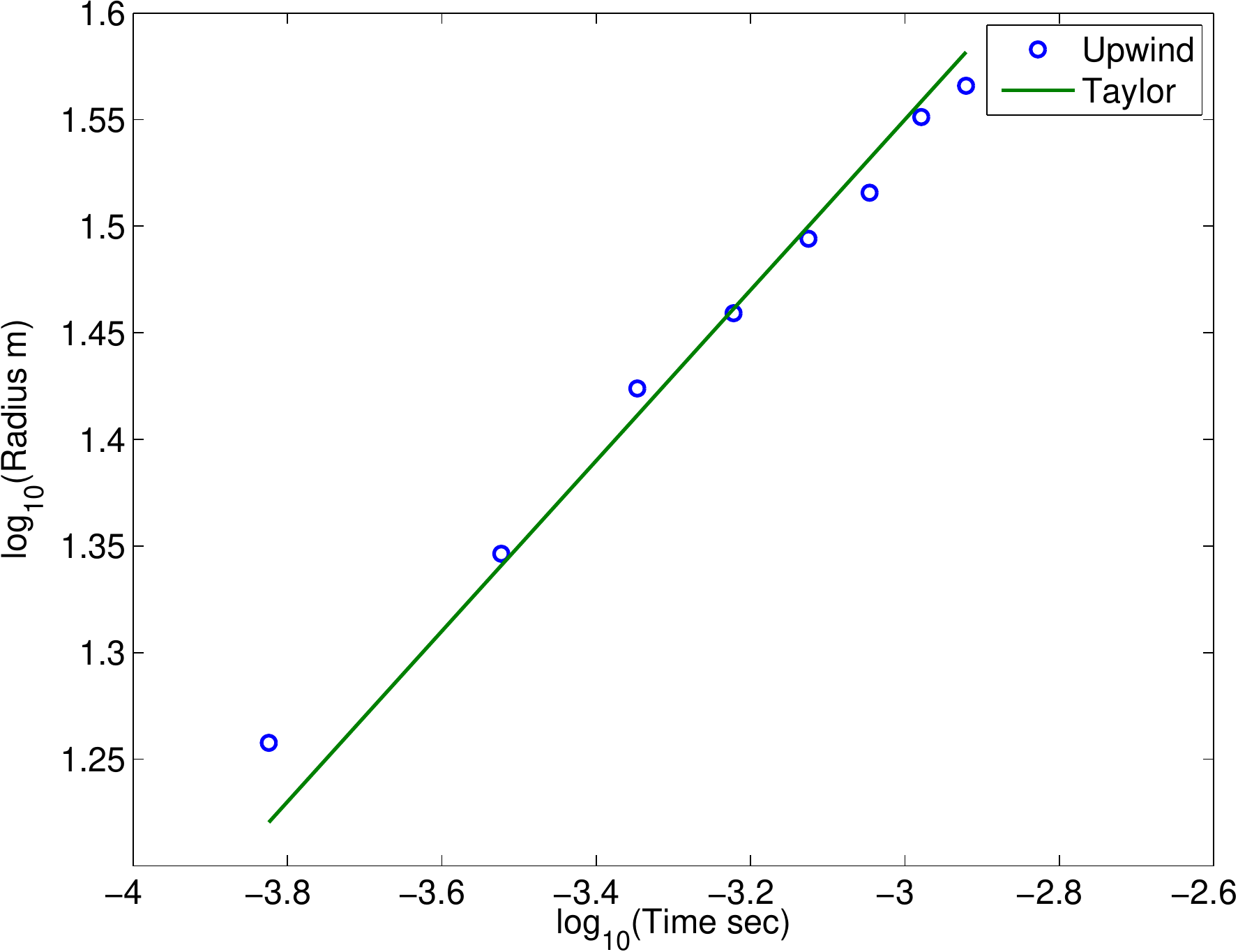} 
\caption{Variation of blast radius with time. The straight line corresponds to the relation $R^{5/2} t^{-1}=$const.}
\label{fig:blastrad}
\end{center}
\end{figure}

The computations are performed in a time accurate manner using 3-stage Runge-Kutta scheme and a constant time step of $6 \times 10^{-8}$ for which the CFL number is less than one. The Frink and Jameson reconstruction schemes lead to loss of positivity of pressure in the first iteration itself. Only the upwind scheme is able to give stable solutions for this problem. In figure~(\ref{fig:blastp}) we compare the radial pressure variation for the first order scheme and the second order scheme with upwind reconstruction at two different times; note that the pressure jump across the shock has a large magnitude in the higher order scheme. It is seen that the second order scheme is also free of oscillations. Moreover, the second order scheme is able to resolve more variations in the solution as compared to first order scheme, especially at the earlier times. This indicates that the limiter is not completely reducing the scheme to first order. A similarity solution for the radius of the blast wave is given by Taylor~\cite{1950RSPSA.201..159T} according to which the quantity $R^{5/2} t^{-1}$ is a constant, where $R$ is the radius of the blast wave. From the numerical solution of the pressure variation obtained with the upwind reconstructed scheme, we determine the blast wave radius and the results are plotted in figure~(\ref{fig:blastrad}) along with the analytical results. Except for the early times, the evolution of the radius is well predicted by the current computations as compared to theoretical results.

\section{Conclusions}

Vertex-centroid finite volume schemes on tetrahedral grids have a simple structure even for second order schemes since they make use of vertex values for reconstruction. Maintaining the positivity of the interpolation weights is important for stability of the scheme. The new averaging procedure given here is shown to be more successful in maintaining positivity of the interpolation weights and has good properties in terms of having O(1) determinants. The simplified reconstruction scheme (upwind scheme) proposed here is shown to be theoretically more accurate as compared to the Frink scheme, and gives comparably similar results on test cases. By writing the scheme as a difference of vertex values and cell-center values, limited versions of the scheme are developed and  shown to be stable in maximum norm for scalar conservation laws. Upwind scheme is seen to be less restrictive in terms of the conditions that the limiter has to satisfy. These schemes give oscillation free solutions when applied to Euler equations. For the blast wave problem which has a very strong shock, only the upwind limited scheme was stable while the other schemes failed due to loss of positivity of density/pressure, demonstrating the good stability properties of the proposed reconstruction scheme.

\bibliographystyle{abbrv}
\bibliography{bibdesk}

\end{document}